\documentclass{amsart}

\input xypic
\usepackage{amsmath,amssymb,amsthm}

\newtheorem{theorem}{Theorem}[section]
\newtheorem{proposition}[theorem]{Proposition}
\newtheorem{lemma}[theorem]{Lemma}
\newtheorem{corollary}[theorem]{Corollary}

\theoremstyle{definition}

\newcommand{\hlgy}[1]{\ensuremath{H_{*}(#1)}}

\newcounter{bean}

\newcommand{\namedright}[3]{\ensuremath{#1\stackrel{#2}
 {\longrightarrow}#3}}
\newcommand{\nameddright}[5]{\ensuremath{#1\stackrel{#2}
 {\longrightarrow}#3\stackrel{#4}{\longrightarrow}#5}} 
\newcommand{\namedddright}[7]{\ensuremath{#1\stackrel{#2}
 {\longrightarrow}#3\stackrel{#4}{\longrightarrow}#5
  \stackrel{#6}{\longrightarrow}#7}} 

\newcommand{\larrow}{\relbar\!\!\relbar\!\!\rightarrow} 
\newcommand{\llarrow}{\relbar\!\!\relbar\!\!\larrow}

\newcommand{\lnameddright}[5]{\ensuremath{#1\stackrel{#2}
 {\larrow}#3\stackrel{#4}{\larrow}#5}} 
 
\newcommand{\llnamedright}[3]{\ensuremath{#1\stackrel{#2}
 {\llarrow}#3}}
\newcommand{\llnameddright}[5]{\ensuremath{#1\stackrel{#2}
 {\llarrow}#3\stackrel{#4}{\llarrow}#5}} 
\newcommand{\llnamedddright}[7]{\ensuremath{#1\stackrel{#2}
 {\llarrow}#3\stackrel{#4}{\llarrow}#5
  \stackrel{#6}{\llarrow}#7}}

\newcommand{\qqed}{\hfill\Box}  
 
\newcommand{\zmodp}{\ensuremath{\mathbb{Z}\mathit{/p}\mathbb{Z}}} 
  
\newcommand{\SW}{\mathcal{W}} 
\newcommand{\Z}{\mathbb{Z}} 
\newcommand{\noqed}{\def\qedsymbol{}}

\begin{document} 

%%% Title 

\title{An elementary construction of Anick's fibration}
\author{Brayton Gray}
\address{Department of Mathematics, Statistics and Computer Science\\
         University of Illinois at Chicago\\
         851 S.~Morgan Street\\
         Chicago, IL, 60607-7045, USA} 
\email{brayton@uic.edu}
\author{Stephen Theriault} 
\address{Department of Mathematical Sciences,
         University of Aberdeen, Aberdeen AB24 3UE, United Kingdom}
\email{s.theriault@maths.abdn.ac.uk}

%\date{}
%\keywords{} 

\begin{abstract}
Cohen, Moore, and Neisendorfer's work on the odd primary homotopy theory of 
spheres and Moore spaces, as well as the first author's work on 
the secondary suspension, predicted the existence of a $p$-local fibration 
\(\nameddright{S^{2n-1}}{}{T}{}{\Omega S^{2n+1}}\) 
whose connecting map is degree $p^{r}$. 
In a long and complex monograph, Anick constructed such a 
fibration for $p\geq 5$ and $r\geq 1$. Using new methods we give a much 
more conceptual construction which is also valid for $p=3$ and $r\geq 1$. 
We go on to establish several properties of the space $T$.  
\end{abstract}

\maketitle

\section{Introduction} 
\label{sec:intro} 

In \cite{CMN1, CMN2, N1} Cohen, Moore, and Neisendorfer proved a landmark 
result concerning the exponent of the homotopy groups of spheres 
localized at an odd prime $p$. When $p\geq 3$ and $r\geq 1$ they constructed
a map 
\(\pi_{n}:\namedright{\Omega^2 S^{2n+1}}{}{S^{2n-1}}\) 
such that the composition with the double suspension 
\[\nameddright{\Omega^2 S^{2n+1}}{\pi_n}{S^{2n-1}}{E^2}{\Omega^2 S^{2n+1}}\] 
is homotopic to the $p^{r}$-power map. The existence of
such a map for $r=1$ was used to show that~$p^{n}$ annihilates the 
$p$-torsion in $\pi_{\ast}(S^{2n+1})=0$.

In \cite{CMN3}, the authors raised the question of whether the map $\pi_{n}$ 
occurs in a fibration sequence
\begin{equation} 
  \label{eqA}
  \tag{A}\namedddright{\Omega^2 S^{2n+1}}{\pi_{n}}{S^{2n-1}}{}{T}  
      {}{\Omega S^{2n+1}}. 
\end{equation}
The first construction of such a fibration was accomplished for $p\geq 5$ 
by Anick~\cite{A} and was the subject of a 270 page book. There has been much 
interest in finding a simpler construction. It is the purpose of 
this paper to give an elementary construction of the space $T$ and 
the fibration~(\ref{eqA}) which is valid for all odd primes. The methods  
are new and have the advantage of being straightforward and accessible 
to nonexperts. It is anticipated that they should be of use for other 
problems as well. A comparison of our methods and Anick's will be 
given once we state our results.  

The question of the existence of a fibration as in~(\ref{eqA}) appeared in another 
context at about the same time. In trying to understand the secondary 
suspension~\cite{C,M}, the first author~\cite{G4,G5} was led to conjecture 
the existence of $(i-1)$-connected spaces $T_i$ which fit into secondary 
$EHP$ sequences
\[\nameddright{T_{2n-1}}{E}{\Omega T_{2n}}{H}{B W_{n}}\] 
\[\nameddright{T_{2n}}{E^{\prime}}{\Omega T_{2n+1}}{H^{\prime}}{BW_{n+1}}\]  
where $BW_{n}$ is the classifying space of the fiber of the double 
suspension constructed in~\cite{G3}. These $EHP$ fibrations should fit 
together in such a way that the resulting spectrum $\{T_i\}$ is
equivalent to the Moore spectrum $S^0\cup_{p^r} e^{1}$. The $T_{i}$'s 
would then give a refinement of the secondary suspension into $2p$ 
stages. The analysis indicated that $T_{2n}$ is homotopy equivalent to 
$S^{2n+1}\{p^{r}\}$, the fiber
of the map of degree $p^{r}$ on $S^{2n+1}$, and that $T_{2n-1}$ would sit 
in the fibration sequence~(\ref{eqA}). 

Our first objective is to construct a secondary Hopf invariant  
\(H:\namedright{\Omega S^{2n+1}\{p^{r}\}}{}{BW_{n}}\) 
for $p\geq 3$. This lets us  
define~$T$ as the homotopy fiber of $H$. It follows easily that $T$  
satisfies the fibration in~(\ref{eqA}) and the second $EHP$ 
fibration. We also show that the space we construct is homotopy 
equivalent to Anick's when $p\geq 5$.   

The $EHP$ viewpoint also predicted that the $T_{i}$'s should have  
a rich structure. They should be homotopy associative and homotopy 
commutative $H$-spaces enjoying a certain universal property. 
Together, these properties would imply that the mod-$p^{r}$ 
homotopy classes of the $T_{i}$'s could be represented 
by multiplicative maps. That is, letting $P^{i}(p^{r})$ be the 
mod-$p^{r}$ Moore space of dimension~$i$, there should be a one-to-one 
correspondence 
\[[P^{i}(p^{r}), T_{j}]\leftrightarrow\{\mbox{homotopy classes of $H$-maps 
     from $T_{i}$ to $T_{j}$}\}.\] 
The properties were easy to establish when $i$ is even~\cite{G4}. 
Subsequent to Anick's work, Anick and the first author~\cite{AG} 
constructed an $H$-space structure on $T$ by showing that, 
for each $n$, there is  
a $(2n-2)$-connected co-$H$ space $G$ with the property that $T$ 
is a retract of $\Omega G$ and $G$ is a retract of $\Sigma T$. 
They also proved a semi-universal property for $T$. The other 
properties were later established by the second author~\cite{T2}. 

Our second objective is to take advantage of our construction of the space $T$  
to give a new, simpler construction of the space $G$, and prove 
all the properties in~\cite{AG} for $p\geq 3$. Collectively, our 
results are as follows.   

\begin{theorem} 
\label{maintheorem} 
   Suppose $p\geq 3$ and $r\geq 1$. Then the following hold: 
   \begin{enumerate}
   \item[(a)] there is an $H$-fibration sequence 
        \[\namedddright{\Omega^2 S^{2n+1}}{\pi_{n}}{S^{2n-1}} 
           {}{T}{}{\Omega S^{2n+1}}\] 
        where the composition 
        \[\nameddright{\Omega^2 S^{2n+1}}{\pi_{n}}{S^{2n-1}}{E^2} 
            {\Omega^2 S^{2n+1}}\] 
        is the $p^{r}$-power map; 
   \item[(b)] there is a fibration sequence 
        \[\namedddright{\Omega G}{h}{T}{}{R}{}{G}\] 
        where $h$ has a right homotopy inverse 
        \(g:\namedright{T}{}{\Omega G}\) 
        so that 
        \[\Omega G\simeq T\times \Omega R\]
        with $R$ a wedge of mod-$p^s$ Moore spaces for $s\geq r$; 
   \item[(c)] the adjoint of $g$, 
        \[\tilde{g}:\namedright{\Sigma T}{}{G},\] 
        has a right homotopy inverse 
        \(f:\namedright{G}{}{\Sigma T}\) 
        and there is a homotopy equivalence 
        \[\Sigma T\simeq G\vee W\]
        where $W$ is a wedge of mod-$p^s$ Moore spaces for $s\geq r$; 
   \item[(d)] there are ``$EHP$ fibrations''
        \[\namedddright{W_n}{P}{T_{2n-1}}{E}{\Omega T_{2n}}{H}{BW_n}\] 
        \[\namedddright{W_{n+1}}{P^{\prime}}{T_{2n}}{E^{\prime}} 
            {\Omega T_{2n+1}}{H^{\prime}}{BW_{n+1}}\] 
        where $T_{2n} =S^{2n+1}\{p^r\}$, $T_{2n-1}=T$, and there  
        is an equivalence of spectra $\{T_i\}\simeq S^0\cup_{p^r} e^1$. 
   \end{enumerate} 
\end{theorem} 

Our methods are simpler and more direct than those of Anick. 
He constructed $T$ as a retract
of a loop space $\Omega D$, where $D$ is an infinite dimensional $CW$-complex
whose bottom two cells are the mod-$p^{r}$ Moore space $P^{2n+1}(p^{r})$ 
and whose other cells come from iteratively attaching certain Moore spaces 
in a delicately prescribed fashion. A great deal of his effort was directed 
towards constructing the attaching maps, and this necessitated the 
introduction of many new techniques. The restriction to primes strictly 
larger than $3$ was due to a heavy reliance on 
differential graded Lie algebras which require that the primes $2$ and $3$ 
be inverted in order for the Lie identities to be satisfied.  
By contrast, we construct the space $T$ directly for all $p\geq 3$ 
without reference to the space $D$ and without reference to differential 
graded Lie algebras. The main ingredient in this new construction is an 
extension theorem (presented as Theorem~\ref{extension}). This allows for  
a straightforward extension of the map 
\(\namedright{\Omega^{2} S^{2n+1}}{}{BW_{n}}\) 
constructed in~\cite{G3} to an $EHP$ map 
\(H:\namedright{\Omega S^{2n+1}\{p\}}{}{BW_{n}}\).  

The new methods may be useful in positively 
resolving a long-standing conjecture  
that the fiber $W_{n}$ of the double suspension 
is a double loop space at odd primes. Including dimension and torsion 
parameters, the space $T_{2np-1}(p)$ gives a candidate for a double  
delooping: potentially $W_{n}\simeq\Omega^{2} T_{2np-1}(p)$. 
Such a homotopy equivalence would have deep implications in homotopy 
theory, one of which being a much better understanding of the differentials 
in the $EHP$ spectral sequence calculating the homotopy groups of spheres.  

This paper is the result of combining separate efforts by the two authors. 
The second author discovered the extension theorem and obtained part~(a) 
of Theorem~\ref{maintheorem} without the $H$-space structure, as well as part~(d). 
The first author later found a different application of the extension theorem 
to obtain a factorization of the map $H$, as well as a further application 
of the extension theorem to obtain parts~(b), (c), and the $H$-space structure.

\section{The extension theorem} 
\label{sec:extension} 

We begin by restating a theorem of the first author~\cite{G3} which identifies 
certain homotopy pullbacks as homotopy pushouts. 
A homotopy fibration 
\(\nameddright{X}{}{Q}{}{A}\) 
has a \emph{trivialization} if there is a homotopy equivalence 
$Q\simeq A\times X$ in which the map 
\(\namedright{Q}{}{A}\) 
becomes the projection 
\(\namedright{A\times X}{\pi_{1}}{A}\).  

\begin{theorem}
   \label{DL}
   Suppose
   \(\nameddright{X}{}{F^{\prime}}{}{E^{\prime}}\)
   is a homotopy fibration and there is a map 
   \(\namedright{A}{}{E^{\prime}}\). 
   Let~$Q$ be the homotopy pullback
   \[\diagram
          Q\rto\dto & F^{\prime}\dto \\
          A\rto & E^{\prime} 
     \enddiagram\] 
   and let $E$ be the homotopy cofiber of 
   \(\namedright{A}{}{E^{\prime}}\). 
   Then the homotopy fibration 
   \(\nameddright{X}{}{Q}{}{A}\) 
   has a trivialization if and only if there is a homotopy pullback
   \[\diagram
        F^{\prime}\rto\dto & F\dto \\
        E^{\prime}\rto & E
     \enddiagram\]
   for some space $F$. Further, if the trivialization exists, there is a 
   homotopy pushout  
   \[\diagram
        Q\simeq A\times X\rto\dto^-{\pi_{2}} & F^{\prime}\dto \\
        X\rto & F
     \enddiagram\]
   where $\pi_{2}$ the projection onto the second factor.
   $\qqed$
\end{theorem} 

There is a special case of Theorem~\ref{DL} in the context of principal 
fibrations which is the key tool used to construct $T$ and 
prove Theorem~\ref{maintheorem}. 
In general, suppose there is a homotopy fibration sequence 
\[\namedddright{\Omega B}{\partial}{F}{}{E}{}{B}.\]  
Then there is a canonical homotopy action 
\(\theta:\namedright{F\times\Omega B}{}{F}\) 
satisfying homotopy commutative diagrams 
\[\diagram 
      \Omega B\times\Omega B\rto^-{\mu}\dto^{\partial\times 1} 
        & \Omega B\dto^{\partial} 
        & & F\times\Omega B\rto^-{\theta}\dto^{\pi_{1}} & F\dto \\ 
      F\times\Omega B\rto^-{\theta} & F 
        & & F\rto & E 
  \enddiagram\] 
where $\mu$ is the loop multiplication and $\pi_{1}$ is the projection. 
Note that both squares are homotopy pullbacks. 
Now suppose there is a homotopy cofibration 
\(\nameddright{A}{b}{E^{\prime}}{}{E}\). 
Define spaces $Q$ and $F^{\prime}$ by the iterated homotopy pullback diagram 
   \[\diagram 
        Q\rto\dto & F^{\prime}\rto\dto & F\dto \\ 
        A\rto^-{b}\dto & E^{\prime}\rto\dto & E\dto \\ 
        B\rdouble & B\rdouble & B. 
     \enddiagram\] 
In particular, the map 
\(\namedright{A}{}{B}\) 
is null homotopic as it factors through the middle row which consists 
of two consecutive maps in a 
homotopy cofibration. So $Q\simeq A\times\Omega B$. This lets us 
apply Theorem~\ref{DL} to see that there is a homotopy pushout 
\[\diagram 
      A\times\Omega B\rto\dto^{\pi_{2}} & F^{\prime}\dto \\ 
      \Omega B\rto & F. 
  \enddiagram\] 
What we wish to do is choose a particular trivialization of $Q$ which 
lets us identify the map 
\(\namedright{A\times\Omega B}{}{F^{\prime}}\). 

The fact that there is some decomposition $Q\simeq A\times\Omega B$ 
implies that we can choose a lift  
\[a:\namedright{A}{}{F^{\prime}}\] 
of $b$. There may be many choices of a lift, but for the moment any choice 
suffices. The definition of $F^{\prime}$ as a homotopy pullback results 
in a homotopy fibration sequence 
\(\namedddright{\Omega B}{}{F^{\prime}}{}{E^{\prime}}{}{B}\).  
This determines a homotopy action 
\(\theta:\namedright{F^{\prime}\times\Omega B}{}{F^{\prime}}\). 
Let $\overline{\theta}$ be the composite 
\[\overline{\theta}:\nameddright{A\times\Omega B}{a\times 1} 
     {F^{\prime}\times\Omega B}{\theta}{F^{\prime}}.\] 

\begin{proposition} 
   \label{DLcase}  
   Let 
   \(\nameddright{F}{}{E}{}{B}\) 
   be a homotopy fibration and suppose there is a homotopy cofibration 
   \(\nameddright{A}{b}{E^{\prime}}{}{E}\). 
  Define the space $F^{\prime}$ and the 
   map $\overline{\theta}$ as above. Then there is a homotopy pushout  
   \[\diagram 
        A\times\Omega B\rto^-{\overline{\theta}}\dto^{\pi_{2}} 
           & F^{\prime}\dto \\ 
        \Omega B\rto & F. 
     \enddiagram\] 
\end{proposition} 

\begin{proof} 
Consider the diagram 
\[\diagram 
     A\times\Omega B\rto^-{a\times 1}\dto^{\pi_{1}} 
         & F^{\prime}\times\Omega B\rto^-{\theta}\dto^{\pi_{1}} 
         & F^{\prime}\dto \\ 
     A\rto^-{a} & F^{\prime}\rto & E^{\prime}. 
  \enddiagram\] 
The right square is a homotopy pullback as it is one of the canonical properties 
of the homotopy action~$\theta$. The left square is a homotopy pullback by the 
naturality of the projection. So the outer rectangle is also a homotopy 
pullback. Observe that the top row of the rectangle is the definition 
of~$\overline{\theta}$ while the bottom row is the given map $b$ 
by the definition of $a$ as a lift. Thus if $Q$ is the homotopy pullback 
\[\diagram 
       Q\rto^-{f}\dto^-{g} & F^{\prime}\dto \\ 
       A\rto & E^{\prime}  
  \enddiagram\] 
then there is a homotopy equivalence 
\(e:\namedright{A\times\Omega B}{}{Q}\) 
such that $g\circ e\sim\pi_{1}$ -- so the homotopy fibration 
\(\nameddright{\Omega B}{}{Q}{}{A}\) \
has been trivialized -- and $g\circ e\sim\overline{\theta}$. Therefore 
Theorem~\ref{DL} implies the existence of the asserted homotopy pushout. 
\end{proof} 

We now state Theorem~\ref{extension}, which uses Proposition~\ref{DLcase} 
to construct an extension under certain conditions. The conditions 
involve exponent information, so we first make two definitions. If $A$ is 
a co-$H$ space, let 
\(\underline{p}^{r}:\namedright{A}{}{A}\) 
be the map of degree~$p^{r}$. If $Z$ is an $H$-space, let 
\(p^{r}:\namedright{Z}{}{Z}\) 
be the $p^{r}$-power map. 

\begin{theorem} 
   \label{extension} 
   Let 
   \[\diagram 
       \Omega B\rto\ddouble & F^{\prime}\rto\dto & E^{\prime}\rto\dto 
          & B\ddouble \\ 
       \Omega B\rto & F\rto & E\rto & B 
     \enddiagram\] 
   be a homotopy fibration diagram and suppose there is a 
   homotopy cofibration 
   \(\nameddright{A}{b}{E^{\prime}}{}{E}\) 
   where $A$ is a suspension. Observe that the map $b$ lifts  
   to $F^{\prime}$; suppose there is a choice of lift   
   \(a:\namedright{A}{}{F^{\prime}}\)
   with the property that
   $\Sigma a\sim t\circ\underline{p}^{r}$ for some map $t$. 
   Suppose there is a map   
   \(f^{\prime}:\namedright{F^{\prime}}{}{Z}\) 
   where $Z$ is a homotopy associative $H$-space whose $p^{r}$-power  
   map is null homtopic. Then there is an extension 
   \[\diagram 
         F^{\prime}\rto^-{f^{\prime}}\dto & Z\ddouble \\ 
         F\rto^-{f} & Z  
     \enddiagram\] 
   for some map $f$. 
\end{theorem} 

Before beginning the proof, we state a Theorem of James~\cite{J} 
and prove two preliminary Lemmas. If $X$ is a space, let 
\(E:\namedright{X}{}{\Omega\Sigma X}\) 
be the suspension.  

\begin{theorem}
   \label{James}
   Let $X$ be a path-connected space and $Z$ be a homotopy
   associative $H$-space. Let
   \(f:\namedright{X}{}{Z}\)
   be a map. Then there is a unique $H$-map
   \(\overline{f}:\namedright{\Omega\Sigma X}{}{Z}\)
   such that $\overline{f}\circ E\sim f$.~$\qqed$
\end{theorem} 

We say that $\overline{f}$ is the \emph{multiplicative extension} of $f$. 

To prepare for Lemmas~\ref{projfactor} and~\ref{expfactor} we establish 
some notation. Let $X$ and $Y$ be spaces. Let 
\(i_{1}:\namedright{X}{}{X\times Y}\) 
and 
\(i_{2}:\namedright{X}{}{X\times Y}\) 
be the inclusions, and let 
\(\pi_{1}:\namedright{X\times Y}{}{X}\) 
and 
\(\pi_{2}:\namedright{X\times Y}{}{Y}\) 
be the projections. It is well known that there is a natural homotopy  
equivalence 
\[e:\namedright{\Sigma X\vee\Sigma Y\vee(\Sigma X\wedge Y)}
    {}{\Sigma (X\times Y)}\]
such that the restrictions of $e$ to $\Sigma X$ and $\Sigma Y$ are 
$\Sigma i_{1}$ and $\Sigma i_{2}$, and $\Sigma\pi_{1}\circ e$ and 
$\Sigma\pi_{2}\circ e$ are homotopic to the pinch maps onto $\Sigma X$ 
and $\Sigma Y$. There may be many choices of such a homotopy equivalence; 
any fixed choice will do. Let $j$ be the restriction   
\[j:\Sigma X\vee(\Sigma X\wedge Y)\hookrightarrow\namedright 
     {\Sigma X\vee\Sigma Y\vee(\Sigma X\wedge Y)}{e}{\Sigma(X\times Y)}.\]    

\begin{lemma} 
   \label{projfactor} 
   Let $Z$ be a homotopy associative $H$-space. Suppose there is a map  
   \(f:\namedright{X\times Y}{}{Z}\) 
   whose multiplicative extension 
   \(\overline{f}:\namedright{\Omega\Sigma(X\times Y)}{}{Z}\) 
   has the property that the composite 
   \[\namedddright{X\vee(X\wedge Y)}{E}{\Omega\Sigma(X\vee(X\wedge Y))} 
         {\Omega j}{\Omega\Sigma(X\times Y)}{\overline{f}}{Z}\]  
   is null homotopic. Then there is a homotopy commutative diagram 
   \[\diagram 
        \Omega\Sigma(X\times Y)\rto^-{\overline{f}}\dto^{\Omega\Sigma\pi_{2}} 
           & Z\ddouble \\ 
        \Omega\Sigma Y\rto^-{\overline{f}_{Y}} & Z 
     \enddiagram\] 
   where $\overline{f}_{Y}$ is the multiplicative extension of 
   \(f_{Y}:\nameddright{Y}{i_{2}}{X\times Y}{f}{Y}\). 
\end{lemma} 

\begin{proof} 
Consider the diagram 
\[\diagram 
     \Omega\Sigma (X\vee Y\vee (X\wedge Y)) 
            \rto^-{\Omega e}\dto^{\Omega\Sigma q} 
        & \Omega\Sigma(X\times Y) 
            \rto^-{\overline{f}}\dto^{\Omega\Sigma\pi_{2}} 
        & Z\ddouble \\ 
     \Omega\Sigma Y\rdouble & \Omega\Sigma Y\rto^-{\overline{f}_{Y}} 
        & Z 
  \enddiagram\] 
where $q$ is the pinch map. By definition of $e$, we have  
$\Sigma\pi_{2}\circ e\sim\Sigma q$, so  the left square homotopy 
commutes. The assertion of the Lemma is that the right square homotopy 
commutes as well. As $e$ is a homotopy equivalence, it is equivalent 
to show that the outer rectangle homotopy commutes. Since all maps are 
multiplicative and $Z$ is homotopy associative, Theorem~\ref{James} 
implies that it is equivalent to show that 
$\overline{f}\circ\Omega e\circ E\sim 
    \overline{f}_{Y}\circ\Omega\Sigma q\circ E$,  
where 
\(E:\namedright{X\vee Y\vee(X\wedge Y)}{} 
     {\Omega\Sigma(X\vee Y\vee(X\wedge Y))}\) 
is the suspension. By hypothesis and the naturality of $E$, the restriction of 
$\overline{f}\circ\Omega e\circ E$ to $X\vee(X\wedge Y)$ is null 
homotopic, so $\overline{f}\circ\Omega e\circ E$ factors as the composite 
\(\nameddright{X\vee Y\vee (X\wedge Y)}{q}{Y}{f_{Y}}{Z}\). 
On the other hand, as $\overline{f}_{Y}$ is the multiplicative extension  
of $f_{Y}$, we have $\overline{f}_{Y}\circ E\sim f_{Y}$. The naturality of 
the suspension therefore implies that 
$\overline{f}_{Y}\circ\Omega\Sigma q\circ E\sim 
    \overline{f}_{Y}\circ E\circ q\sim f_{Y}\circ q$. 
Hence 
$\overline{f}\circ\Omega e\circ E\sim f_{Y}\circ q\sim   
   \overline{f}_{Y}\circ\Omega\Sigma q\circ E$,  
as required. 
\end{proof} 

\begin{lemma} 
   \label{expfactor} 
   Let 
   \(\namedddright{\Omega B}{\partial}{F}{}{E}{}{B}\) 
   be a homotopy fibration sequence and let 
   \(\theta:\namedright{F\times\Omega B}{}{F}\) 
   be the associated homotopy action. Suppose $A$ is a suspension  
   and there is a map 
   \(a:\namedright{A}{}{F}\)  
   such that $\Sigma a\sim t\circ\underline{p}^{r}$ for some map $t$. 
   Let $\overline{\theta}$ be the composite 
   \[\overline{\theta}:\nameddright{A\times\Omega B}{a\times 1} 
        {F\times\Omega B}{\theta}{F}.\] 
   Suppose there is a map 
   \(f:\namedright{F}{}{Z}\) 
   where $Z$ is a homotopy associative $H$-space whose $p^{r}$-power 
   map is null homotopic. Then there is a homotopy commutative diagram 
   \[\diagram 
        A\times\Omega B\rto^-{\overline{\theta}}\dto^{\pi_{2}} 
          & F\dto^{f} \\ 
        \Omega B\rto^-{f\circ\partial} & Z. 
     \enddiagram\] 
\end{lemma} 

\begin{proof} 
First suspend and look at $\Sigma(a\times 1)$. Consider the diagram 
\[\diagram 
    \Sigma A\vee(\Sigma A\wedge\Omega B) 
           \rrto^-{\Sigma\underline{p}^{r}\vee 
               (\Sigma\underline{p}^{r}\wedge 1)}\ddouble  
      & & \Sigma A\vee(\Sigma A\wedge\Omega B) 
           \dto^{t\vee (t\wedge 1)} \\ 
    \Sigma A\vee(\Sigma A\wedge\Omega B)  
           \rrto^-{\Sigma a\vee(\Sigma a\wedge 1)}\dto^{j} 
      & & \Sigma F\vee(\Sigma F\wedge\Omega B)\dto^{j} \\
    \Sigma(A\times\Omega B)\rrto^-{\Sigma(a\times 1)} 
      & & \Sigma(F\times\Omega B).   
  \enddiagram\] 
The top square homotopy commutes by the hypothesis that 
$\Sigma a\sim t\circ\underline{p}^{r}$. The bottom square homotopy 
commutes by the naturality of the map $j$. 

Looping this diagram and using the naturality of the suspension, we 
obtain a homotopy commutative diagram 
\begin{equation} 
  \label{bigprdgrm} 
  \diagram 
     A\vee(A\wedge\Omega B) 
            \rrto^-{\underline{p}^{r}\vee(\underline{p}^{r}\wedge 1)}\dto^{E} 
        & & A\vee (A\wedge\Omega B)\dto^{E} \\  
     \Omega\Sigma (A\vee(A\wedge\Omega B))  
      \rrto^-{\Omega\Sigma (\underline{p}^{r}\vee(\underline{p}^{r}\wedge 1))} 
      \dto^{\Omega j}   
        & & \Omega\Sigma (A\vee(A\wedge\Omega B)) 
           \dto^{\Omega (j\circ(t\vee (t\wedge 1)))} \\ 
     \Omega\Sigma (A\times\Omega B) 
           \rrto^-{\Omega\Sigma (a\times 1)}  
        & & \Omega\Sigma (F\times\Omega B).  
  \enddiagram 
\end{equation}  
Let 
$\phi=\Omega(j\circ(t\vee(t\wedge 1)))\circ E\circ 
    (\underline{p}^{r}\vee(\underline{p}^{r}\wedge 1))$ 
be the upper direction around~(\ref{bigprdgrm}), and let 
$\varphi=\Omega\Sigma(a\times 1)\circ\Omega j\circ E$ be the lower 
direction around~(\ref{bigprdgrm}). So $\phi\sim\varphi$. 
Now compose to $Z$ as follows. By hypothesis, $Z$ is 
homotopy associative, so by Theorem~\ref{James} the identity map on $Z$ 
extends to an $H$-map 
\(r:\namedright{\Omega\Sigma Z}{}{Z}\) 
such that $r\circ E\sim 1$. Define $\gamma$ by the composite 
\[\gamma:\llnamedddright{\Omega\Sigma(F\times\Omega B)}
   {\Omega\Sigma\theta_{k-1}}{\Omega\Sigma F}{\Omega\Sigma f}
   {\Omega\Sigma Z}{r}{Z}.\] 
Observe that $\gamma\circ\phi$ is an element of 
\[N=[A\vee(A\wedge\Omega B), Z]\]  
which is divisible by $p^{r}$. Here, the group structure on $N$ 
is determined by 
$A\vee(A\wedge\Omega B)$ being a suspension. 
As $Z$ is homotopy associative, this group structure on $N$ is 
equivalent to the one determined by the $H$-structure on $Z$. 
By hypothesis, the $p^{r}$-power map on $Z$ is null homotopic, and so $N$ 
has exponent $p^{r}$. Thus $\gamma\circ\phi$ is null homotopic, 
and so $\gamma\circ\varphi$ is null homotopic. 
 
Now we set up to apply Lemma~\ref{projfactor}. Consider the diagram 
\[\diagram 
     A\times\Omega B\rto^-{a\times 1}\dto^{E} 
       & F\times\Omega B\rto^-{\theta}\dto^{E} 
       & F\rto^-{f}\dto^{E} & Z\dto^{E}\drdouble & \\  
     \Omega\Sigma(A\times\Omega B)\rto^-{\Omega\Sigma(a\times 1)} 
       & \Omega\Sigma(F\times\Omega B)\rto^-{\Omega\Sigma\theta} 
       & \Omega\Sigma F\rto^-{\Omega\Sigma f} 
       & \Omega\Sigma Z\rto^-{r} & Z. 
  \enddiagram\] 
The three squares homotopy commute by the naturality of $E$. The 
right triangle homotopy commutes by the definition of $r$. So the 
entire diagram homotopy commutes. Let $g=f\circ\theta\circ(a\times 1)$ 
be the composite along the top row, and let 
$\overline{g}=r\circ\Omega\Sigma f\circ\Omega\Sigma\theta\circ 
   \Omega\Sigma(a\times 1)$     
be the composite along the bottom row. Observe that $\overline{g}$ is an 
$H$-map as it is the composite of $H$-maps, and $\overline{g}\circ E\sim g$. 
Thus $\overline{g}$ is the multiplicative extension of $g$. Further, 
by their definitions, $\overline{g}=\gamma\circ\Omega\Sigma(s\times 1)$ and  
$\varphi=\Omega\Sigma(s\times 1)\circ\Omega j\circ E$, so the null 
homotopy for $\gamma\circ\varphi$ implies that the composite 
\[\lnameddright{A\vee(A\wedge\Omega B)}{\Omega j\circ E} 
     {\Omega\Sigma(A\times\Omega B)}{\overline{g}}{Z}\] 
is null homotopic. Therefore, by Lemma~\ref{projfactor}, there is a      
homotopy commutative square 
\[\diagram 
      \Omega\Sigma(A\times\Omega B) 
           \rto^-{\overline{g}}\dto^{\Omega\Sigma\pi_{2}} 
         & Z\ddouble \\ 
      \Omega\Sigma\Omega B\rto^-{\overline{h}} & Z 
  \enddiagram\]   
where $\overline{h}$ is the multiplicative extension of 
\[h:\nameddright{\Omega B}{i_{2}}{A\times\Omega B}{g}{Z}.\]  

Finally, the previous square and the naturality of the suspension give 
a homotopy commutative diagram 
\[\diagram 
      A\times\Omega B\rto^-{E}\dto^{\pi_{2}} 
         & \Omega\Sigma(A\times\Omega B) 
             \rto^-{\overline{g}}\dto^{\Omega\Sigma\pi_{2}} 
         & Z\ddouble \\ 
      B\rto^-{E} & \Omega\Sigma\Omega B\rto^-{\overline{h}} & Z. 
  \enddiagram\] 
Since $\overline{g}$ and $\overline{h}$ are the multiplicative extensions 
of $g$ and $h$ respectively, we have $\overline{g}\circ E\sim g$ and 
$\overline{h}\circ E\sim h$, and so the homotopy commutativity of the 
diagram implies that $g\sim h\circ\pi_{2}$. Now we untangle definitions. 
By their definitions, $g=f\circ\theta\circ(a\times 1)$ and 
$\overline{\theta}=\theta\circ(a\times 1)$. So $g\sim f\circ\overline{\theta}$. 
By their definitions, $h=g\circ i_{2}$ and $g=f\circ\theta\circ(a\times 1)$.  
The definition of $\theta$ as a homotopy action implies that  
$\theta\circ(a\times 1)\circ i_{2}\sim\partial$. 
Thus $h\sim f\circ\partial$. Hence  
$f\circ\overline{\theta}\sim f\circ\partial\circ\pi_{2}$, 
precisely as asserted by the Lemma.  
\end{proof}  
 
\noindent\textit{Proof of Theorem~\ref{extension}:} 
The given diagram of principal fibrations and homotopy 
cofibration let us apply Proposition~\ref{DLcase} to obtain a 
homotopy pushout 
\[\diagram 
     A\times\Omega B\rto^-{\overline{\theta}}\dto^{\pi_{2}} 
       & F^{\prime}\dto \\ 
    \Omega B\rto & F. 
  \enddiagram\] 
Since $A$ is a suspension, the lift 
\(\namedright{A}{a}{F^{\prime}}\) 
has the property that $\Sigma a\sim t\circ\underline{p}^{r}$, 
and $Z$ is a homotopy associative $H$-space whose $p^{r}$-power 
map is null homotopic, we can apply  
Lemma~\ref{expfactor} to the homotopy fibration sequence 
\(\namedddright{\Omega B}{}{F^{\prime}}{}{E^{\prime}}{}{B}\) 
and the given map 
\(\namedright{F^{\prime}}{f^{\prime}}{Z}\) 
in order to obtain a homotopy commutative diagram   
\[\diagram 
     A\times\Omega B\rto^-{\overline{\theta}}\dto^{\pi_{2}} 
        & F^{\prime}\dto^{f^{\prime}} \\ 
    \Omega B\rto & Z. 
  \enddiagram\] 
Therefore there is a pushout map 
\(f:\namedright{F}{}{Z}\) 
with the property that the composite 
\(\nameddright{F^{\prime}}{}{F}{f}{Z}\) 
is homotopic to $f^{\prime}$, as required. 
$\qqed$

\section{The construction of the space $T$} 
\label{sec:Tconstruction} 

The purpose of this section is to construct the spaces
$T$ and produce several fibration sequences. 
We begin our discussion with the Moore space
\[P^k(p^r)=S^{k-1}\cup_{p^r}e^k\]
which we will abbreviate as $P^k$. Let us fix some notation
by defining a diagram of fibration sequences
induced by the lower right hand corner
\begin{equation} 
  \label{dgrmB}\tag{B} 
  \diagram 
      \Omega^{2} S^{2n+1}\rto^-{\partial}\dto & E\rto^-{\pi}\dto^{\sigma} 
         & F\rto\dto & \Omega S^{2n+1}\dto \\ 
      \ast\rto\dto & P^{2n+1}\rdouble\dto & P^{2n+1}\rto\dto & \ast\dto \\ 
      \Omega S^{2n+1}\rto & S^{2n+1}\{p^{r}\}\rto & S^{2n+1}\rto^-{p^{r}} 
          & S^{2n+1}. 
    \enddiagram 
\end{equation}
The spaces $E$ and $F$ were first introduced in \cite{CMN2,CMN1}. It
is easy to see that
\[H^i(F)=\begin{cases}
          \Z&i=2kn\\
          0&i\neq 2kn
         \end{cases}\]  
and the fibration connecting map 
\(\namedright{\Omega S^{2n+1}}{}{F}\) 
is divisible by $p^{r}$ in each nonzero degree in integral cohomology.  
In their work \cite{CMN2}, the authors introduced certain maps 
\(x_i:\namedright{P^{2ni-1}}{}{\Omega F}\) 
whose adjoints
\(\tilde{x}_i:\namedright{P^{2ni}}{}{F}\) 
induce epimorphisms in integral cohomology. For each
$i>1$, $x_i$ is a relative Samelson product, so the composition 
\[\nameddright{P^{2ni}}{\tilde{x}_{i}}{F}{}{}{P^{2n+1}}\] 
is an iterated Whitehead product. Since $S^{2n+1}\{p^r\}$ is
an $H$ space, these classes lift to $E$, giving diagrams  
\[\diagram 
     P^{2ni}\rto^-{y_{i}}\dto^{\tilde{x}_{i}} & E\dto^{\sigma} \\ 
     F\rto & P^{2n+1} 
  \enddiagram\] 
for some maps $y_{i}$. In particular, 
\(\pi y_i -\tilde{x}_i:\namedright{P^{2ni}}{}{F}\) 
composes trivially to $P^{2n+1}$ and so factors through $\Omega S^{2n+1}$. 
Thus the induced homomorphism in $2ni$ dimensional
cohomology is divisible by~$p^r$ and hence trivial.
Consequently, we obtain the following. 

\begin{lemma} 
   \label{PEFcohlgy} 
   The composite 
   \(\nameddright{H^{2ni}(F)}{\pi^{\ast}}{H^{2ni}(E)}{y_i^{\ast}}
       {H^{2ni}(P^{2ni})}\) 
   is an epimorphism for \linebreak each $i>1$.~$\qqed$ 
\end{lemma} 

We require one more lemma to apply the results of
Section~\ref{sec:extension}. 

\begin{lemma} 
   \label{AMcone}
   Suppose $X$ is 2-connected and $M$ is
   either a sphere or a Moore space. Let 
   \(f:\namedright{\Sigma M}{}{X}\) 
   be given. Define $A$ by the cofibration sequence 
   \[\namedddright{M}{p^s}{M}{i}{A}{j}{\Sigma M}.\] 
   Suppose there is a commutative diagram 
   \[\diagram 
        \Sigma A\rto^-{\Sigma j}\dto^{x} & \Sigma^{2} M\ddouble \\ 
        X\cup_{f} C\Sigma M\rto^-{\rho} & \Sigma^{2} M 
     \enddiagram\] 
   for some map $x$, where $\rho$ is the quotient map. Then there is 
   a commutative diagram 
   \[\diagram 
        \Sigma M\rto^-{\Sigma i}\dto^{x^{\prime}} & \Sigma A\dto^{x} \\ 
        X\rto & X\cup_{f} C\Sigma M 
     \enddiagram\] 
   for some map $x^{\prime}$, and $f$ is homotopic to $p^{s}\cdot x^{\prime}$. 
\end{lemma} 

\begin{proof}
Consider the standard map from a cofibration
sequence to a fibration sequence defined by the right
hand square 
\[\diagram 
     \Sigma M\rto^-{p^{s}}\dto^{x^{\prime\prime}} 
        & \Sigma M\rto^-{\Sigma i}\dto^{x^{\prime}} 
        & \Sigma A\rto^-{\Sigma j}\dto^{x} & \Sigma^{2} M\ddouble\\ 
     J(\Sigma M)\rto & J(X,\Sigma M)\rto 
        & X\cup_{f} C\Sigma M\rto^-{\rho} & \Sigma^{2} M. 
  \enddiagram\] 
Here $J(\Sigma M)$ is the James construction and $J(X,SM)$ is the fiber 
of $\rho$~\cite{G2}. The map $x^{\prime\prime}$ is the adjoint to the
identity for an appropriate choice of $x^{\prime}$. Suppose $\Sigma M$ has
dimension $k$. Since $X$ is $2$-connected, the $k+1$
skeleton of $J(X,\Sigma M)$ is $X$ and $x^{\prime}$ factors through $X$
up to homotopy. Since~$x^{\prime\prime}$ factors through $\Sigma M$ as well
we have a homotopy commutative square 
\[\diagram 
     \Sigma M\rto^-{p^{s}}\ddouble & \Sigma M\dto^{x^{\prime}} \\ 
     \Sigma M\rto^-{f} & X 
  \enddiagram\] 
which proves the Lemma. 
\end{proof} 

We apply these results as follows. Let $F_{(i)}$ be the
$2ni$ skeleton of $F$, so
\[F_{(i)} =F_{(i-1)}\cup_{\gamma_i} e^{2ni}\] 
where $\gamma_{i}$ is the attaching map. 
Now combining \ref{PEFcohlgy} and \ref{AMcone} with $M=S^{2ni-2}$,
$X=F_{(i-1)}$, $f=\gamma_i$, $s=r$, and $x=\pi y_i$ we obtain the following. 

\begin{corollary} 
   \label{Fdivbyp}
   For each $i>1$ we have a homotopy commutative diagram
   \[\diagram 
         S^{2ni-1}\rto^-{p^{r}}\ddouble & S^{2ni-1}\dto^{\delta_{i}} \\ 
         S^{2ni-1}\rto^-{\gamma_{i}} & F_{(i-1)} 
     \enddiagram\] 
   where $\delta_i$ satisfies a homotopy commtutative diagram 
   \[\diagram 
        S^{2ni-1}\rto\dto^{\delta_{i}} & P^{2ni}\rto & E\dto^{\pi} \\ 
        F_{(i-1)}\rrto & & F. 
     \enddiagram\] 
   $\qqed$ 
\end{corollary} 

We now set up to apply Theorem~\ref{extension}. Define the space $E_{(i)}$ 
as the homotopy pullback 
\[\diagram 
     E_{(i)}\rto\dto & E\dto^-{\pi} \\ 
     F_{(i)}\rto & F. 
  \enddiagram\] 
Observe that there is a homotopy pullback diagram 
\[\diagram 
     \Omega^{2} S^{2n+1}\rto\ddouble & E_{(i-1})\rto\dto 
         & F_{(i-1)}\rto\dto & \Omega S^{2n+1}\ddouble \\ 
     \Omega^{2} S^{2n+1}\rto & E_{(i)}\rto & F_{(i)}\rto 
         & \Omega S^{2n+1}. 
  \enddiagram\] 
By the definition of $F_{(i)}$ there is a homotopy cofibration 
\(\nameddright{S^{2ni-1}}{\gamma_{i}}{F_{(i-1)}}{}{F_{(i)}}\). 
Thus $\gamma_{i}$ lifts to $E_{(i-1)}$. A lift can be chosen 
which is divisible by $p^{r}$. Specifically, by Corollary~\ref{Fdivbyp},  
$\gamma_{i}\sim\delta_{i}\circ p^{r}$. Moreover, 
\(\namedright{S^{2ni-1}}{\delta_{i}}{F_{(i-1)}}\)  
composed to $F$ factors through 
\(\namedright{E}{\pi}{F}\). 
Thus there is a pullback map 
\(\bar{y}_{i}:\namedright{S^{2ni-1}}{}{E_{(i-1)}}\) 
such that the composite 
\(\nameddright{S^{2ni-1}}{\bar{y}_{i}}{E_{(i-1)}}{}{F_{(i-1)}}\) 
is homotopic to~$\delta_{i}$.  
Hence $a=\bar{y}_{i}\circ p^{r}$ is a lift of $\gamma_{i}$.   
Theorem~\ref{extension} now immediately implies the following. 

\begin{theorem} 
   \label{Eiextension}
   If $Z$ is a homotopy associative $H$ space whose $p^{r}$-power map 
   is null homotopic, then for $i>1$  any map 
   \(\namedright{E_{(i-1)}}{\phi}{Z}\) 
   extends to a map 
   \(\overline{\phi}:\namedright{E_{(i)}}{}{Z}\). 
   $\qqed$ 
\end{theorem} 

In~\cite{G3}, a classifying space $BW_{n}$ of the fiber of the double suspension 
was constructed, along with a fibration sequence 
\[\nameddright{S^{2n-1}}{E^{2}}{\Omega^{2} S^{2n+1}}{\nu}{BW_{n}}.\]  

\begin{corollary} 
   \label{EtoBWn}
   There is a map 
   \(\nu^E:\namedright{E}{}{BW_n}\) 
   such that the composition
   \[\nameddright{\Omega^2 S^{2n+1}}{\partial}{E}{\nu^E}{BW_n}\] 
   is homotopic to $\nu$.
\end{corollary} 

\begin{proof}
Since $F_{(1)} =S^{2n}$, we have the fibration 
\[\namedddright{\Omega^2 S^{2n+1}}{}{E_{(1)}}{}{S^{2n}}{}{\Omega S^{2n+1}}.\]

This fibration was analyzed in~\cite{G3} and it
was shown that $E_{(1)}\simeq S^{4n-1}\times BW_n$ in such a way
that the composition 
\[\nameddright{\Omega^2 S^{2n+1}}{\partial}{S^{4n-1}\times BW_n} 
    {\pi_2}{BW_n}\] 
is homotopic to $\nu$. It was also shown that for $p\geq 5$
$BW_n$ is a homotopy associative $H$ space. The $H$ space 
structure on $BW_n$ was shown to be homotopy
associative for $p=3$  and that the $p^{th}$-power map on $BW_n$ is 
null homotopic in~\cite{T5}. Thus for $i>1$ we can apply 
Theorem~\ref{Eiextension} to construct maps
\(\nu_i:\namedright{E_{(i)}}{}{BW_n}\) 
by induction such that 
$\nu_i\partial_i \sim \nu$. Since $E=\cup E_{(i)}$, we define 
\(\nu^E:\namedright{E}{}{BW_n}\) 
by $\nu^E\mid E_i=\nu_i$.
\end{proof} 

\begin{theorem} 
   \label{Fsplit}
There is a diagram of fibrations 
\[\diagram 
    S^{2n-1}\rto\dto^{i} & \Omega^{2} S^{2n+1}\rto^-{\nu}\dto^{\partial} 
      & BW_{n}\ddouble \\ 
    R_{0}\rto\dto & E\rto^-{\nu^{E}}\dto & BW_{n} \\ 
    F\rdouble & F & 
  \enddiagram\] 
with $i$ null homotopic and so $\Omega F\simeq S^{2n-1}\times \Omega R_0$.
\end{theorem} 

\begin{proof}
The space $R_0$ is defined as the fiber of $\nu^E$. Since the
fibration 
\[\nameddright{\Omega^2 S^{2n+1}}{\partial}{E}{}{F}\] 
is induced by a map to $\Omega S^{2n+1}$ which induces an
isomorphism in $H_{2n}(\ )$, the map $\Omega F\xrightarrow{}S^{2n-1}$ induces
an isomorphism in $H_{2n-1}(\ )$ and hence has a right homotopy inverse.
\end{proof}

It is worth noting at this point that the space
$\Omega R_0$ is split in \cite{CMN1}; there is a homotopy decomposition 
\[\Omega R_0\simeq \prod_{i\geq 1} S^{2np^i-1}\{p^{r+1}\} 
    \times\Omega P (n,r)\]
where $P(n,r)$ is a complicated wedge of mod-$p^r$ Moore
spaces. The fact that the product on the right is a loop space and is
mapped to $\Omega F$ by a loop map is not
obvious from their analysis. The structure of $R_0$ is
rather simple. 

\begin{proposition} 
   \label{Rcohlgy} 
   We have 
   \[H^m(R_0)=\begin{cases}
       \Z/{p^ri} &\text{if}\ m=2ni\\
       0&\text{otherwise.}
       \end{cases}\] 
Furthermore, there is a choice of generators
$e_i\in H^{2mi}(X)$ such that $e_ie_j=p^r\binom{i+j}{i} e_{i+j}$.
\end{proposition} 

\begin{proof}
Apply the Serre spectral sequence to the fibration 
\(\nameddright{S^{2n-1}}{}{R_{0}}{}{F}\) 
in Theorem~\ref{Fsplit}. 
\end{proof} 

We now construct the space $T$ in Theorem~\ref{maintheorem} 
and prove the existence of the fibrations in parts~(a) and~(d), 
leaving the $H$-structure to the next section. 
By Diagram~(\ref{dgrmB}) there is a fibration sequence 
\(\namedddright{\Omega S^{2n+1}\{p^{r}\}}{\tau}{E}{\sigma}{P^{2n+1}}{} 
    {S^{2n+1}\{p^{r}\}}\). 
Define $H$ by the composition 
\[H:\nameddright{\Omega S^{2n+1}\{p^{r}\}}{\tau}{E}{\nu^{E}}{BW_{n}}.\] 
Note that $H$ can be regarded as a secondary Hopf invariant. 
Define $T$ as the homotopy fiber of~$H$. Then Theorem~\ref{Fsplit} 
implies the following. 

\begin{theorem} 
   \label{Tfib}
   There is a diagram of fibrations
   \[\diagram 
         T\rto\dto & \Omega S^{2n+1}\{p^{r}\}\rto^-{H}\dto^{\tau} 
            & BW_{n}\ddouble \\ 
         R_{0}\rto\dto & E\rto^-{\nu^{E}}\dto^{\sigma} & BW_{n} \\ 
         P^{2n+1}\rdouble & P^{2n+1}. & 
     \enddiagram\] 
   $\qqed$ 
\end{theorem} 

The connecting maps for the vertical fibrations in 
Theorem~\ref{Tfib} immediately give the following. 

\begin{corollary} 
   \label{loopPtoT}   
   There is a homotopy commutative diagram 
   \[\diagram 
        \Omega P^{2n+1}\rdouble\dto & \Omega P^{2n+1}\dto \\ 
        T\rto & \Omega S^{2n+1}\{p^{r}\} 
     \enddiagram\] 
   where the right map is the loop of the inclusion of the bottom 
   Moore space.~$\qqed$ 
\end{corollary} 

Continuing the diagram in~(\ref{dgrmB}), we have 
\[\diagram 
     \Omega^{2} S^{2n+1}\rto^-{\rho}\ddouble 
        & \Omega S^{2n+1}\{p^{r}\}\dto^{\tau} \\ 
     \Omega^{2} S^{2n+1}\rto^-{\partial} & E. 
  \enddiagram\] 
Observe that $H\rho\sim\nu^{E}\tau\rho\sim\nu^{E}\partial\sim\nu$. 
Theorem~\ref{Tfib} therefore implies the following. 

\begin{theorem} 
   \label{CMNfib}
   There is a diagram of fibrations 
   \[\diagram 
        \Omega^{2} S^{2n+1}\rdouble\dto^{\pi_{n}} 
           & \Omega^{2} S^{2n+1}\dto^{p^{r}} & \\ 
        S^{2n-1}\rto^-{E^{2}}\dto & \Omega^{2} S^{2n+1}\rto^-{\nu}\dto^{\rho} 
           & BW_{n}\ddouble \\ 
        T\rto\dto & \Omega S^{2n+1}\{p^{r}\}\rto^-{H}\dto 
           & BW_{n} \\ 
        \Omega S^{2n+1}\rdouble & \Omega S^{2n+1}. & 
     \enddiagram\] 
   $\qqed$ 
\end{theorem}

In particular, the top square in Theorem~\ref{CMNfib} is Cohen, 
Moore, and Neisendorfer's factorization of the $p^{r}$-power map 
on $\Omega^{2} S^{2n+1}$. Since $\pi_{n}$ has degree $p^{r}$, we 
have the following corollary. 

\begin{corollary} 
   \label{CMNsquare}
   There is a homotopy commutative diagram
   \[\diagram 
        S^{2n-1}\rto^-{p^{r}}\dto & S^{2n-1}\dto \\ 
        \Omega^{2} S^{2n+1}\rto^-{p^{r}}\urto^-{\pi_{n}} 
           & \Omega^{2} S^{2n+1} 
     \enddiagram\] 
   for each $r\geq 1$.~$\qqed$ 
\end{corollary}

\section{The construction of $G$ and the $H$-space structure on $T$} 
\label{sec:G} 

In this section we construct an $H$-space structure on $T$. In
fact we do more than that. We construct a corresponding
co-$H$ space~$G$ in the sense of \cite{G7}; i.e., we construct a 
$(2n-2)$-connected space $G$ and maps
\begin{align*}
   & f:\namedright{G}{}{\Sigma T} \\
   & g:\namedright{T}{}{\Omega G} \\
   & h:\namedright{\Omega G}{}{T}  
\end{align*}
such that the compositions
\begin{align*} 
   & \nameddright{G}{f}{\Sigma T}{\tilde{g}}{G} \\ 
   & \nameddright{T}{g}{\Omega G}{h}{T} 
\end{align*}
are homotopic to the identity, where $\tilde{g}$ is the adjoint of
$g$. We go on to derive several interesting results from this structure. 

We will write $T^m$ for the $m$-skeleton of $T$. We will also reintroduce 
the torsion parameter for Moore spaces as we will need   
to consider mod-$p^s$ Moore spaces $P^{m}(p^{s})$ 
for $s\neq r$. The space $G$ will be filtered by subcomplexes $G_k$ 
which will be constructed inductively starting with $G_{-1}=\ast$. We will
construct a map 
\[\alpha_{k}:\namedright{P^{2np^k}(p^{r+k})}{}{G_{k-1}}\] 
and define $G_k$ as the mapping cone of $\alpha_k$.

The induction proceeds through 14 steps for each $k$, and we collect some 
information outside of the induction first. 

\begin{proposition} 
   \label{Tcohlgy}
   As an algebra, $H^{\ast}(T; \Z/p)$ is generated by classes
   $u$ of dimension $2n-1$ and $v_i$ of dimension $2np^i$ for each $i\geq 0$
   subject to the relations $v_i^p=0$ and $u^2=0$. For each $i$ define
   \[u_i=uv_0^{p-1}v_1^{p-1}\dots v_{i-1}^{p-1}\]
   Then $\beta^{(r+i)} u_i=v_i$. As a vector space $\tilde{H}^{\ast}(T;\Z/p)$ 
   is generated by classes $v(m)$ of dimension $2mn$ and $u(m)$ of dimension
   $2mn-1$ for each $m\geq 1$ where
   \begin{align*}
      v(m)&=v_s^{e_s}\dots v_t^{e_t}=\beta^{(r+s)}u(s)\\
      u(m)&=u_sv_s^{e_s}v_{s+1}^{e_{s+1}}\dots v_t^{e_t}
   \end{align*}
   and $m=\sum^t_{i=s}e_ip^i$, $0\leq e_i<p$, $e_s\neq 0$.
\end{proposition} 

\begin{proof}
We apply the Serre spectral sequence for the
cohomology of the fibration 
\[\nameddright{S^{2n-1}}{}{T}{}{\Omega S^{2n+1}}\]
Using $\Z/p$ coefficients we see that
\[H^{\ast}(T;\Z/p) \cong H^{\ast} \left(S^{2n-1}; \Z/p\right)\otimes
     H^{\ast}\left(\Omega S^{2n+1};\Z/p\right)\]
as algebras. Using integer coefficients we see that $\nu(m)$
is the reduction of a class of order $p^{r+s}$ so $v(m)=\beta^{(r+s)}u(s)\neq 0$.
We define $v_s=\beta^{(r+s)}u_s$.
\end{proof}

Note that dually the homology of $T$ has a very simple description. There 
is a Hopf algebra isomorphism 
\[\hlgy{T}\cong\Lambda(\bar{u})\otimes\zmodp [\bar{v}]\] 
where $\bar{u}$ and $\bar{v}$ are dual to $u$ and $v$ respectively, 
and the dual Bocksteins are determined by 
$\beta^{(r+i)}\bar{v}^{p^{i}}=\bar{u}\bar{v}^{p^{i}-1}$ for $i\geq 0$.    

Anick \cite{A} introduced the  notation $\SW^b_a$ for the class of all  
spaces that are locally finite wedges of mod-$p^s$ Moore spaces 
for $a\leq s\leq b$. Note that any simply connected Moore space is a suspension, 
so any simply connected space in $\SW^{b}_{a}$ is a suspension. Recall that the smash 
of two Moore spaces is homotopy equivalent to a wedge of Moore spaces: if $s\leq t$ 
then there is a homotopy equivalence 
\[P^{m}(p^{s})\wedge P^{n}(p^{t})\simeq P^{m+n}(p^{s})\wedge P^{m+n-1}(p^{s}).\] 
In particular, $\SW^{b}_{a}$ is closed under smash products. Recall also that 
any retract of a wedge of Moore spaces is homotopy equivalent to a wedge of 
Moore spaces, so $\SW^{b}_{a}$ is closed under retracts. 

\begin{lemma} 
   \label{wedgelemma}
   Suppose $W\in \SW^b_a$ is simply connected and 
   \(f:\namedright{P^k(p^t)}{}{W}\) 
   is divisible by $p^b$. 

   \textup{(a)}\hspace*{0.5em} Write $W=W_1\vee W_2$ with 
   $W_1\in \SW^{b-1}_a$ and $W_{2}\in\SW^{b}_{b}$. Then $f$
   factors through $W_2$ up to homotopy.

   \textup{(b)}\hspace*{0.5em} Suppose in addition that $W_2$ is 
   $(d-1)$ connected and $k<pd$. Then $f\sim\ast$.
\end{lemma} 

\begin{proof} 
Since $W$ is a wedge, there is a homotopy equivalence 
$\Omega W=\Omega W_2\times \Omega(W_1\rtimes \Omega W_2)$ (see, for example,~\cite{G1}). 
Since $W_{1},W_{2}\in\SW^{b}_{a}$, both spaces are suspensions, and we can write 
$W_{1}=\Sigma\overline{W}_{1}$ and $W_{2}=\Sigma\overline{W}_{2}$. Since $W_{1}$ is 
a suspension, we have 
$W_{1}\rtimes\Omega W_{2}\simeq W_{1}\vee (W_{1}\wedge\Omega W_{2})$. 
For the right wedge summand, the James splitting of $\Sigma\Omega\Sigma X$ as 
$\bigvee\Sigma X^{(i)}$ gives 
\[W_{1}\wedge\Omega W_{2}\simeq\Sigma\overline{W}_{1}\wedge\Omega\Sigma\overline{W}_{2} 
    \simeq\overline{W}_{1}\wedge\left(\bigvee\Sigma \overline{W}_{2}^{(i)}\right).\]  
Combining, we have 
\[W_{1}\rtimes\Omega W_{2}\simeq 
    W_{1}\vee\left(W_{1}\wedge\left(\bigvee W_{2}^{(i)}\right)\right).\] 
In particular, since $\SW^{b}_{a}$ is closed under smash products, we have  
$W_{1}\rtimes\Omega W_{2}\in\SW^{b}_{a}$. Applying the 
Hilton-Milnor theorem therefore implies that 
$\Omega(W_1\rtimes\Omega W_2)\simeq\prod_i\Omega P^{n_i}(p^{s_{i}})$ 
with $a\leq s\leq b-1$. 

By~\cite{N3}, the $p^{r+1}$-power map on $\Omega^{2} P^{m}(p^{r})$ is null homotopic 
for any $r\geq 1$ and $m\geq 3$. Thus $P^{m}(p^{r})$ admits no nontrivial maps which 
are divisible by $p^{r+1}$. In our case, this implies that $\prod_{i}\Omega P^{n_{i}}(p^{s_{i}})$ 
admits no nontrivial maps which are divisible by $p^{b}$. Thus the adjoint of $f$, 
which is divisible by $p^{b}$, is trivial on $\Omega (W_{1}\rtimes\Omega W_{2})$ and 
so factors through the inclusion 
\(\namedright{\Omega W_{2}}{}{\Omega W}\). 
Hence, adjointing, $f$ factors through the inclusion 
\(\namedright{W_{2}}{}{W}\), 
proving part~(a).  

For part~(b), since $W_{2}\in\SW^{b}_{b}$ and $W_{2}$ is $(d-1)$-connected, 
the Hilton-Milnor theorem implies that 
$\Omega W_2=\prod\Omega P^{n_i}(p^b)$ where $n_i > d$. By~\cite{CMN1,N3}, 
$P^{2m+1}(p^{r})$ admits no nontrivial maps which are divisible by $p^{r}$ from 
a $CW$-complex of dimension $t<2mp$, and $P^{2m}(p^{b}))$ admits no nontrivial maps 
which are divisible by $p^{r}$ from a $CW$-complex of dimension $t<2(2m-1)p$. In our case, 
the $CW$-complex is $P^{k}(p^{t})$, the domain of $f$, and the target Moore spaces 
are the $P^{n_{i}}(p^{b})$ in the decomposition of $\Omega W_{2}$. Since $n_{i}>d$ 
for each $i$, the hypothesis $k<pd$ guarantees that the component of $f$ 
on $P^{n_{i}}(p^{b})$, being divisible by $p^{b}$, is null homotopic. Hence $f$ 
is null homotopic. 
\end{proof} 

\begin{theorem} 
   \label{Ginduct}
   For each $k\geq 0$ there are spaces $G_k$ and $W_k\in \SW_r^{r+k-1}$
   satisfying the following conditions:
   \begin{enumerate}
   \item[(a)]
   $\Sigma T^{2np^k-2}\simeq G_{k-1}\vee W_k$; 
   \item[(b)]
   there are maps 
   \(g_k:\namedright{T^{2np^k-2}}{}{\Omega G_{k-1}}\) 
   and
   \(h_{k-1}:\namedright{\Omega G_{k-1}}{}{T}\) 
   such that $h_{k-1}g_k$ is homotopic to the inclusion of
   $T^{2np^k-2}$ into $T$; 
   \item[(c)]
   there is a homotopy commutative diagram of cofibration
   sequences which defines $G_k$ 
   \[\diagram 
       P^{2np^{k}}(p^{r+k})\rto^-{m_{k}}\ddouble 
         & \Sigma T^{2np^{k}-2}\rto\dto^{\tilde{g}_{k}} 
         & \Sigma T^{2np^{k}}\dto^{g^{\prime}_{k}} \\ 
       P^{2np^{k}}(p^{r+k})\rto^-{\alpha_{k}} & G_{k-1}\rto & G_{k} 
     \enddiagram\] 
   where $\tilde{g}_k$ is the adjoint of $g_k$; 
   \item[(d)]
   there is a map 
   \(e:\namedright{P^{2np^k}(p^{r+k-1})\vee P^{2np^k+1}(p^{r+k-1})} 
       {}{\Sigma T^{2np^k}}\) 
   which induces an epimorphism in mod-$p$ cohomology; 
   \item[(e)]
   the map 
   \(m_k:\namedright{P^{2np^k}(p^{r+k})}{}{\Sigma T^{2np^k-2}}\) 
   is divisible by $p^{r+k-1}$; 
   \item[(f)]
   there is a map 
   \(\varphi_k:\namedright{G_k}{}{S^{2n+1}\{p^r\}}\) 
   extending $\varphi_{k-1}$; 
   \item[(g)]
   $\Sigma G_k\in \SW^{r+k}_r$; 
   \item[(h)]
   there is a homotopy commutative diagram of fibration sequences 
   \[\diagram 
       \Omega G_{k}\rto^-{h_{k}} & T\rto\dto & R_{k}\rto\dto 
          & G_{k}\ddouble & \\ 
       & \Omega S^{2n+1}\{p^{r}\}\rto\dto^{H} & E_{k}\rto\dto^{\nu_{k}} 
          & G_{k}\rto & S^{2n+1}\{p^{r}\} \\ 
       & BW_{n}\rdouble & BW_{n}; & 
     \enddiagram\] 
   \item[(i)]
   $\Sigma^{2}\Omega G_{k-1}\in \SW^{r+k-1}_r$; 
   \item[(j)]
   the equivalence in \emph{(a)} extends to an
   equivalence $\Sigma T^{2np^k}\simeq G_k\vee W_k$; 
   \item[(k)]
   $\Sigma^2 T^{2np^k}\in \SW^{r+k}_r$; 
   \item[(l)]
   $G_k\wedge T^{2np^k}\in\SW^{r+k}_r$; 
   \item[(m)]
   $\Sigma T^{2np^k}\wedge T^{2np^k}\in \SW^{r+k}_r$; 
   \item[(n)]
   there is a map 
   \(\mu_k:\namedright{T^{2np^k}\times T}{}{T}\) 
   which is the inclusion on the first axis and the identity on the second.
   Furthermore there is a homotopy commutative square 
   \[\diagram 
        T^{2np^{k}}\times T\rto^-{\mu_{k}}\dto & T\dto \\ 
        \Omega S^{2n+1}\times\Omega S^{2n+1}\rto & \Omega S^{2n+1}.  
     \enddiagram\] 
   \end{enumerate}
\end{theorem} 

\begin{proof}
With $G_{-1}=\ast$ and $G_0=P^{2n+1}$ these statements are
all immediate for $k=0$ with 
\linebreak  
\(\varphi_0:\namedright{P^{2n+1}}{}{S^{2n+1}\{p^r\}}\) 
the inclusion, $E_0=E$ from Theorem~\ref{Fsplit}, $\nu_0=\nu^E$, 
\(\mu_0:\namedright{P^{2n}\times T}{}{T}\) 
obtained from the action of $\Omega P^{2n+1}$ on $T$ defined by the
fibration in Theorem~\ref{Fsplit}. We now supposed that (a)--(n) are all
valid with $k-1$ in the place of $k$ and we proceed to prove
them for $k$.\noqed 
\end{proof} 

\begin{proof}[Proof of {\normalfont (a)}]
We will construct a map
\[f_m:\namedright{P^{2mn+1}(p^{r+s})}{}{\Sigma T^{2np^k-2}}\]
which induces a monomorphism in mod-$p$ homology
for each $m$ satisfying $p^{k-1}< m< p^k$, where $s=\nu_p(m)$. We then
assemble these into a map 
\[\namedright{\Sigma T^{2np^{k-1}}\vee\left( 
   \bigvee^{p^k-1}_{m=p^{k-1}+1}P^{2mn+1}(p^{r+s})\right)} 
   {}{\Sigma T^{2np^k-2}}\]
which induces an isomorphism in mod-$p$ homology. By
applying (j) in the case $k-1$ we are done.

To construct the maps $f_m$ we appeal to (n) in the case $k-1$
and iterate this to produce a diagram with $p$ factors 
\[\diagram 
     T^{2np^{k-1}}\times\dots\times T^{2np^{k-1}}\rto\dto & T^{2np^k}\dto \\ 
     J(S^{2n})_{p^{k-1}}\times\dots\times J(S^{2n})_{p^{k-1}}\rto 
        & J(S^{2n})_{p^k} 
   \enddiagram\] 
where $J(S^{2n})_{j}$ is the $2nj$ skeleton of $\Omega S^{2n+1}$. 
Since $p^{k-1}< m < p^k$,we can write $m=a_sp^s+\dots+a_{k-1}p^{k-1}$ 
with $a_s >0$ and $a_{k-1} >0$.
Write $l=a_sp^s+\dots+a_{k-2}p^{k-2}$ so that $m=l + a_{k-1}p^{k-1}$ and
further restrict the above diagram to one with $a_{k-1}+1$ factors 
\[\diagram 
     T^{2nl}\times T^{2np^{k-1}}\times\dots\times T^{2np^{k-1}} 
        \rto^-{\overline{\mu}}\dto & T^{2nm}\dto \\
     J(S^{2n})_{l}\times J(S^{2n})_{p^{k-1}}\times\dots\times 
        J(S^{2n})_{p^{k-1}}\rto & J(S^{2n})_{m}. 
  \enddiagram\] 
By applying the maps in this diagram to a generator of
$H^{2mn}(J(S^{2n})_{m}; \Z/p)$ we see that
\[(\overline{\mu})^{\ast}\left(v(m)\right) = v(l)\otimes
      v_{k-1}\otimes\dots\otimes v_{k-1}.\]
Now $v(m)=\beta^{(r+s)}u(m)$ and $v(l)=\beta^{(r+s)} u(l)$, so
\[v(l)\otimes v_{k-1}\otimes\dots\otimes
    v_{k-1}=\beta^{(r+s)}\left(u(l)\otimes v_{k-1}\otimes\dots\otimes
    v_{k-1}\right).\] 
Applying (k) and (l) in case $k-1$ we see that
\[\Sigma\left(T^{2nl}\times T^{2np^{k-1}}\times\dots\times T^{2np^{k-1}}\right) 
     \in\SW_r^{r+k-1}.\]
Now given any space $W\in \SW_r^{r+k-1}$ and any class $\xi \in H^i(W;\Z/p)$
with $\beta^{(j)}\xi\neq~0$, there is a map
\[f_{\xi}:\namedright{P^{i+1}(p^j)}{}{W}\] 
with $f_{\xi}^{\ast}$ an epimorphism. Thus for each $m$ satisfying 
$p^{k-1}<m<p^k$ we may choose such a map corresponding to $\xi=u(l)\otimes
v_{k-1}\otimes\dots\otimes v_k$. The composition 
\[\nameddright{P^m(p^{r+s})}{f_{\xi}} 
    {\Sigma\left(T^{2nl}\times T^{2np^{k-1}}\times\dots\times
    T^{2np^{k-1}}\right)}{\Sigma\overline{\mu}}{\Sigma T^{2mn}}\]
therefore gives the desired map $f_m$.\noqed
\end{proof} 

\begin{proof}[Proof of {\normalfont (b)}]
From part (a) we obtain a map 
\(\namedright{T^{2np^k-2}}{}{\Omega G_{k-1}}\) 
which induces an isomorphism in~$\pi_{2n-1}(\ )$. The composition 
\[\nameddright{T^{2np^k-2}}{}{\Omega G_{k-1}}{h_{k-1}}{T}\] 
factors through $T^{2np^k-2}$ and provides a self map of $T^{2np^k-2}$
which induces an isomorphism on~$\pi_{2n-1}(\ )$. Calculations with
cup products and Bocksteins show that this map is a
homotopy equivalence, so composing with the inverse
provides a possibly different map 
\[g_k:\namedright{T^{2np^k-2}}{}{\Omega G_{k-1}}\]
such that $h_{k-1}g_k$ is homotopic to the inclusion.\noqed 
\end{proof} 

\begin{proof}[Proof of {\normalfont (c)}]
Using the map $g_k$ from (b) we construct
a commutative diagram where the bottom row is the 
fibration sequence from (h) in case $k-1$ and the
middle row is a cofibration sequence 
\[\diagram 
    & & T^{2np^k}/T^{2np^k-2}\rdouble\dto & P^{2np^{k}}(p^{r+k})\dto \\ 
    T^{2np^{k}-2}\rto\dto^{g_{k}} & T\rto\ddouble & T/T^{2np^{k}-2}\rto\dto 
       & \Sigma T^{2np^{k}-2}\dto^{\tilde{g}_{k}} \\ 
    \Omega G_{k-1}\rto^-{h_{k-1}} & T\rto & R_{k-1}\rto & G_{k-1}.  
  \enddiagram\] 
Define 
\(\alpha_k:\namedright{P^{2np^{k}}(p^{r+k})}{}{G_{k-1}}\) 
as the vertical composition on the right. Define $g^{\prime}_{k}$ 
by the diagram of cofibration sequences 
\[\diagram 
     P^{2np^k}(p^{r+k})\rto^-{m_k}\ddouble 
       & \Sigma T^{2np^k-2}\rto\dto^{\tilde{g}_{k}} 
       & \Sigma T^{2np^k}\dto^{g^{\prime}_{k}} \\
    P^{2np^k}(p^{r+k})\rto^-{\alpha_{k}} & G_{k-1}\rto & G_{k}.  
  \enddiagram\] 
\noqed
\end{proof} 

\begin{proof}[Proof of {\normalfont (d)}]
As in part (a), we consider the diagram:
\[\diagram 
     T^{2np^{k-1}}\times\dots\times T^{2np^{k-1}}\rto^-{\tilde{\mu}}\dto  
       & T^{2np^k}\dto \\
     J(S^{2n})_{p^{k-1}}\times\dots\times J(S^{2n})_{p^{k-1}}\rto 
       & J(S^{2n})_{p^k} 
  \enddiagram\]
with $p$ factors on the left. This is
defined by iterated application of part (n) in case $k-1$.
Clearly 
\[\tilde{\mu}^{\ast}(v_k)=v_{k-1}\otimes\dots\otimes
     v_{k-1}=\beta^{(r+k-1)}(uv_1^{p-1}\dots
     v^{p-1}_{k-2}\otimes v_{k-1}\otimes\dots\otimes v_{k-1}).\]
As before there is a map 
\[q:\namedright{P^{2np^k+1}(p^{r+k-1})}{} 
    {\Sigma\left(T^{2np^{k-1}}\times\dots\times T^{2np^{k-1}}\right)}\] 
such that $(\Sigma(\tilde\mu) q)^{\ast}$ is an epimorphism in $\Z/p$ cohomology
obtained by applying (k) and (m) in case~$k-1$. Similarly,
$v_ k=\beta^{(r+k)} u_k$ and
\[(\tilde{\mu})^{\ast} u_k 
    =\sum_{p\ \text{terms}} v_{k-1}\otimes\dots\otimes v_{k-1}\otimes
    u_{k-1}\otimes v_{k-1}\dots \otimes v_{k-1}.\]
In particular, the map
\[\namedright{T^{2np^{k-1}-1}\times T^{2np^{k-1}}\times\dots
     \times T^{2np^{k-1}}}{\tilde{\mu}^{\prime}}{T^{2np^k}}\] 
has the property that
\begin{align*}
   \left(\tilde{\mu}^{\prime}\right)^{\ast}(u_k)&=u_{k-1}\otimes
   v_{k-1}\otimes\dots\otimes v_{k-1}\\
   &= \beta^{(r+k-1)}\left(u_{k-1}\otimes u_{k-1}\otimes
   v_k\otimes\dots\otimes v_k\right).
\end{align*}
It follows, as before, that there is a map
\[r:\namedright{P^{2np^k}(p^{r+k-1})}{}
    {\Sigma\left(T^{2np^{k-1}-1}\times T^{2np^{k-1}}\times\dots\times
    T^{2np^{k-1}}\right)}\]
such that $(\Sigma(\tilde{\mu}^{\prime})r)^{\ast}$ is an epimorphism in $\Z/p$
cohomology. We construct $e$ as the wedge sum 
\[e=(\Sigma\tilde{\mu}^{\prime})r\vee (\Sigma\tilde{\mu})q: 
     \namedright{P^{2np^k}(p^{r+k-1})\vee P^{2np^k+1}(p^{r+k-1})} 
      {}{\Sigma T^{2np^k}}.\] 
\noqed
\end{proof} 

\begin{proof}[Proof of {\normalfont (e)}]
We apply Lemma~\ref{AMcone} with $x=e$, $s=r+k-1$, 
$X=\Sigma T^{2np^k-2}$, $M= P^{2np^k-1}(p^{r+k})$, and 
\(f=m_k:\namedright{\Sigma M}{}{X}\). 
In this case
$A=P^{2np^k-1}(p^{r+k-1})\vee P^{2np^k}(p^{r+k-1})$, 
which is the cofiber of $p^{r+k-1}$ on $M=P^{2np^k-1}\left(p^{r+k}\right)$. 
It follows that $m_k$ is divisible by~$p^{r+k-1}$.\noqed
\end{proof} 

\begin{proof}[Proof of {\normalfont (f)}]
To show that there is an extension of $\varphi_{k-1}$ to $\varphi_k$, 
\[\diagram 
    P^{2np^{k}}(p^{r+k})\rto^-{\alpha_{k}} 
      & G_{k-1}\rto\dto^{\varphi_{k-1}} 
      & G_{k}=G_{k-1}\cup_{\alpha_{k}} CP^{2np^{k}}(p^{r+k}) 
          \dlto_{\varphi_{k}} \\ 
    & S^{2n+1}\{p^{r}\} & 
  \enddiagram\] 
it suffices to show that $\alpha_k$ is divisible by $p^r$. This
holds by (e) since $r+k-1\geq r$.\noqed
\end{proof} 

\begin{proof}[Proof of {\normalfont (g)}]
By part (g) for $k-1$, there is a homotopy equivalence   
$SG_{k-1}\simeq\bigvee\limits^{k-1}_{i=0}P^{2np^i+2}\left(p^{r+i}\right)$. 
Also, by definition, 
$SG_k=SG_{k-1}\cup_{S\alpha_k}CP^{2np^k}\left(p^{r+k}\right)$. By part (e)
$\alpha_k=\tilde{\alpha}\circ\left(p^{r+k-1}\iota\right)$ so
$S\alpha_k=\tilde{\alpha}\circ\left(p^{r+k-1}\iota\right)\sim
\left(p^{r+k-1}\iota\right)\circ\tilde{\alpha}$. However
\[
 p^{r+k-1}\iota\colon\bigvee^{k-1}_{i=0}P^{2np^i+2}\left(p^{r+i}\right)\xrightarrow{}\bigvee^{k-1}_{i=0}P^{2np^i+2}\left(p^{r+i}\right)
\]
is null homotopic since the order of the identity map on a mod-$p^{r}$ 
Moore space is $p^{r}$.\noqed
\end{proof} 

\begin{proof}[Proof of {\normalfont (h)}]
As we have constructed 
\(\varphi_k:\namedright{G_k}{}{S^{2n+1}\{p^r\}}\) 
in part (f), we have a pullback diagram of principal fibrations 
\[\diagram 
     \Omega S^{2n+1}\{p^{r}\}\rdouble\dto & \Omega S^{2n+1}\{p^{r}\}\dto \\ 
     E_{k-1}\rto\dto & E_{k}\dto \\ 
     G_{k-1}\rto\dto^{\varphi_{k-1}} & G_{k}\dto^{\varphi_{k}} \\ 
     S^{2n+1}\{p^{r}\}\rdouble & S^{2n+1}\{p^{r}\}. 
  \enddiagram\] 
We wish to apply Theorem~\ref{extension} to extend 
\(\nu_{k-1}:\namedright{E_{k-1}}{}{BW_n}\) 
to $E_k$. It suffices to show that there is a lifting $\delta$ of $\alpha_k$, 
\[\diagram 
    & E_{k-1}\dto \\ 
    P^{2np^{k}}(p^{r+k})\rto^-{\alpha_{k}}\urto^-{\delta} & G_{k-1} 
  \enddiagram\] 
which is divisible by $p$. Since $\alpha_k\sim p^{r+k-1}\tilde{\alpha}$, 
it follows that $p^r\tilde{\alpha}$ lifts to a map 
\(\delta^{\prime}:\namedright{P^{2np^k}(p^{r+k})}{}{E_{k-1}}\) 
with $p^{k-1}\delta^{\prime}=\delta$ a lifting
of $\alpha_k$. Thus as long as $k>1$ we can construct $\delta$ with the
requisite property. When $k=1$, we appeal to~\cite{CMN1} where it
is shown that $\alpha_1=p\delta_1$ with 
\(\delta_1:\namedright{P^{2np}(p^{r+1})}{}{P^{2n+1}(p^r)}\) 
lifting to $E_0$.\noqed
\end{proof} 

\begin{proof}[Proof of {\normalfont (i)}] 
By part (j) in case $k-1$, 
$\Sigma^2\Omega G_{k-1}$ is a retract of $\Sigma^2\Omega\Sigma T^{2np^{k-1}}$. 
The latter space splits since the loop space can be
approximated by the James construction~\cite{J}, giving 
\[\Sigma^2\Omega\Sigma T^{2np^{k-1}}\simeq 
   \Sigma^2\left(\bigvee_{i\geq 1} (T^{2np^{k-1}})^{(i)}\right)\]
which is in $\SW_r^{r+k-1}$ by (k) and (l) in case $k-1$. Since $\SW_r^{r+k-1}$
is closed under retracts we are done.\noqed
\end{proof} 

\begin{proof}[Proof of {\normalfont (j)}]
By part (a), we have $\Sigma T^{2np^k-2}\simeq G_{k-1}\vee W_k$
and by (e), we have
\[\Sigma T^{2np^k}\simeq\left(\Sigma T^{2np^k-2}\right)\cup_{m_k} 
   CP^{2np^k}(p^{r+k})\]
with $m_k$ divisible by $p^{r+k-1}$. It suffices to show that the map
\[m_k:\namedright{P^{2np^k}(p^{r+k})}{} 
    {\Sigma T^{2np^k-2}\simeq G_{k-1}\vee W_k}\]
factors though $G_{k-1}$. To this end, observe that there is a homotopy 
decomposition 
\[\Omega\left(G_{k-1}\vee W_k\right)\simeq 
   \Omega G_{k-1}\times\Omega\left(W_k\rtimes\Omega G_{k-1}\right).\]
We will show that any map
\(\namedright{P^{2np^k}(p^{r+k})}{}{W_k\rtimes\Omega G_{k-1}}\) 
which is divisible by $p^{r+k-1}$ is null homotopic. Since $W_k$ is 
$(4n-1)$-connected, the Moore spaces in $W_{k}$ are double suspensions, so 
$W_k\rtimes \Omega G_{k-1}\in \SW_r^{r+k-1}$. In fact,  
$W_k\rtimes \Omega G_{k-1}\simeq W_1\vee W_2$ with
$W_1\in \SW_r^{r+k-2}$ and $W_2$ a retract of
\[\bigvee^{p-1}_{r=2}P^{2np^{k-1}+1}(p^{r+k-1})\rtimes \Omega G_{k-1}\]
which is $4np^{k-1}-1$ connected. The result follows from
Lemma~\ref{wedgelemma}.\noqed
\end{proof} 

\begin{proof}[Proof of {\normalfont (k)}]
This follows immediately from (g) and (j).\noqed
\end{proof} 

\begin{proof}[Proof of {\normalfont (l)}]
This follows from 3 steps based on an
analysis which first appeared in~\cite{T1}. 

\noindent\textbf{Step 1:}  
$G_k\wedge T^{2np^{k-1}}\in \SW_r^{r+k-1}$. 

Consider the cofibration sequence 
\[\llnameddright{P^{2np^k}(p^{r+k})\wedge T^{2np^{k-1}}} 
    {\alpha_k\wedge 1}{G_{k-1}\wedge T^{2np^{k-1}}}{} 
    {G_k\wedge T^{2np^{k-1}}}.\] 
We have $P^{2np^k}(p^{r+k})\wedge T^{2np^{k-1}}\in \SW_r^{r+k-1}$ 
and $\alpha_k\wedge 1$ is divisible by $p^{r+k-1}$. Consequently, 
$\alpha_k\wedge1\sim \ast$ and so there is a homotopy decomposition 
\[G_k\wedge T^{2np^{k-1}}\simeq (G_{k-1}\wedge T^{2np^{k-1}})\vee
    (P^{2np^k+1}(p^{r+k})\wedge T^{2np^{k-1}})\] 
which is in $\SW_r^{r+k-1}$ by (k) in case $k-1$. 

\noindent\textbf{Step 2:}  
$G_{k-1}\wedge T^{2np^k}\in \SW_r^{r+k-1}$. 

By (j) in case $k-1$, $G_{k-1}\wedge T^{2np^k}$ 
is a retract of $\Sigma T^{2np^{k-1}}\wedge T^{2np^k}$. But
\[\Sigma T^{2np^{k-1}}\wedge T^{2np^k}\simeq
    T^{2np^{k-1}}\wedge\Sigma T^{2np^{k}}\simeq
    T^{2np^{k-1}}\wedge\left(G_k\vee W_k\right)\]
by (j). By Step~1 and (k) in case $k-1$, the latter space is in
$\SW_r^{r+k-1}$. Since $\SW_r^{r+k-1}$ is closed under retracts, we 
therefore have  
$G_{k-1}\wedge T^{2np^k}\in\SW_r^{r+k-1}$. 

\noindent\textbf{Step 3:} 
$G_k\wedge T^{2np^k}\in \SW_r^{r+k}$. 

Consider here the cofibration sequence
\[\llnameddright{P^{2np^k}(p^{r+k})\wedge T^{2np^k}} 
    {\alpha_k\wedge 1}{G_{k-1}\wedge T^{2np^k}}{}{G_k\wedge T^{2np^k}}.\] 
The first space is in $\SW_r^{r+k}$ by (k) and the second is in $\SW_r^{r+k-1}$
by Step~2. In fact, 
$G_{k-1}\wedge T^{2np^k}\simeq 
   (P^{2np^{k-1}+1}(p^{r+k-1})\wedge T^{2np^k})\vee W^{\prime}$ 
with $W^{\prime}\in \SW_r^{r+k-2}$. Here, the projection onto the first 
factor is $\rho_{k-1}\wedge 1$, where $\rho_{k-1}$ is obtained by 
collapsing $G_{k-2}$ to a point.  Applying Lemma~\ref{wedgelemma}(b), 
we see that if $\alpha_k\wedge 1$ is nontrivial, so is the composition 
\[\llnameddright{P^{2np^k}(p^{r+k})\wedge T^{2np^k}}{\alpha_k\wedge 1} 
   {G_{k-1} \wedge T^{2np^k}}{\rho_{k-1}\wedge 1} 
   {P^{2np^{k-1}}(p^{r+k-1})\wedge T^{2np^k}}.\] 
We will show that this composition is null homotopic. 
Let $\delta = \rho_{k-1}\alpha_k$, which is 
divisible by $p^{r+k-1}$ because $\alpha_{k}$ is. According to \cite{N1}, 
the $p^{r+k-1}$-power map on $S^{2np^{k-1}+1}\{p^{r+k-1}\}$ is null 
homotopic. Therefore the composition
\[\nameddright{P^{2np^k}(p^{r+k})}{\delta}{P^{2np^{k-1}+1} (p^{r+k-1})} 
     {}{S^{2np^{k-1}+1}\{p^{r+k-1}\}}\] 
is null homotopic. It follows that the composition 
\[\nameddright{P^{2np^k}(p^{r+k})}{\delta}{P^{2np^{k-1}+1}(p^{r+k-1})} 
   {\rho}{S^{2np^{k-1}+1}}\] 
is null homotopic. Since the map 
\[\llnamedright{P^{2np^{k-1}+1}(p^{r+k})\wedge T^{2np^k}}{\rho\wedge 1}
     {S^{2np^{k-1}+1}\wedge T^{2np^k}}\] 
has a left homotopy inverse, the map 
\[\llnamedright{P^{2np^k}(p^{r+k})\wedge T^{2np^k}}{\delta\wedge 1} 
    {P^{2np^{k-1}+1}(p^{r+k})\wedge T^{2np^k}}\] 
is null homotopic. Since $\alpha_{k}\wedge 1$ is the composition 
\[\llnameddright{P^{2np^k}(p^{r+k})\wedge T^{2np^k}}{\alpha\wedge 1} 
     {G_{k-1}\wedge T^{2np^k}}{\delta\wedge 1} 
     {P^{2np^{k-1}+1}(p^{r+k})\wedge T^{2np^k}},\]
it is null homotopic as well. Consequently, there is a homotopy 
decomposition 
\[G_k\wedge T^{2np^k}\simeq (G_{k-1}\wedge T^{2np^k})\vee 
     (P^{2np^k+1}(p^{r+k})\wedge T^{2np^k}).\] 
Both terms on the right are in $\SW_r^{r+k}$ by (k) and Step~2.\noqed  
\end{proof} 

\begin{proof}[Proof of {\normalfont (m)}]
By (j), 
$\Sigma T^{2np^k}\wedge T^{2np^k}\simeq (G_k\vee W_k) \wedge T^{2np^k}$. 
By (l), $G_{k}\wedge T^{2np^{k}}\in \SW_r^{r+k}$, and as~$W_{k}$ is a 
wedge of Moore spaces which are at least $(4n-1)$-connected, it is a 
double suspension, so by~(k) we have $W_{k}\wedge T^{2np^{k}}\in\SW_{r}^{r+k}$. 
Thus $\Sigma T^{2np^k}\wedge T^{2np^k}\in\SW_{r}^{r+k}$.\noqed
\end{proof} 

\begin{proof}[Proof of {\normalfont (n)}]
Since the composite
\(\namedddright{R_k}{}{E_k}{}{G_k}{}{S^{2n+1}\{p^r\}}\longrightarrow S^{2n+1}\) 
is null homotopic by (h), there is a commutative diagram of
principal fibrations:
\[\diagram 
      \Omega G_{k}\rto\dto^{h_{k}} & \Omega S^{2n+1}\ddouble \\ 
      T\rto\dto & \Omega S^{2n+1}\dto \\ 
      R_{k}\rto\dto & PS^{2n+1}\dto \\ 
      G_{k}\rto & S^{2n+1} 
  \enddiagram\] 
where $PS^{2n+1}$ is the path space on $S^{2n+1}$. 
Consequently the actions are compatible 
\[\diagram 
     \Omega G_k\times T\rto\dto 
        & \Omega S^{2n+1}\times\Omega S^{2n+1}\dto \\
     T\rto & \Omega S^{2n+1}. 
  \enddiagram\] 
Using (j) we construct a map 
\(g_k:\namedright{T^{2np^k}}{}{\Omega G_k}\) 
such that the composition
\[\nameddright{T^{2np^k}}{g_k}{\Omega G_k}{h_k}{T}\] 
is homotopic to the inclusion as in (b). This gives a
homotopy commutative diagram 
\[\diagram 
     T^{2np^k}\times T\rto^-{g_{k}\times 1}\dto^{\mu_{k}} 
       & \Omega G_k\times T\dto^{a} \\
     T\rdouble & T. 
  \enddiagram\] 
Combining the preceeding two diagrams gives the result and completes
the induction.~$\qqed$\noqed 
\end{proof} 

We now consider the limiting case. Write $G=\bigcup G_k$, $R=\bigcup R_k$ 
and $E_{\infty}=\bigcup E_k$. 

\begin{theorem} 
   \label{Glimit}
   There is a diagram of fibration sequences 
   \[\diagram 
        \Omega G\rto^-{h} & T\rto^-{i}\dto^{E} & R\rto\dto 
          & G\ddouble & \\ 
        & \Omega S^{2n+1}\{p^{r}\}\rto\dto^{H} & E_{\infty}\rto\dto 
          & G\rto^-{\varphi} & S^{2n+1}\{p^{r}\} \\ 
        & BW_{n}\rdouble & BW_{n} & & 
     \enddiagram\]  
   and there are maps 
   \(\tilde{g}:\namedright{T}{}{\Omega G}\) 
   and 
   \(f:\namedright{G}{}{\Sigma T}\) 
   such that the composites 
   \[\nameddright{G}{f}{\Sigma T}{g}{G}\] 
   \[\nameddright{T}{\tilde{g}}{\Omega G}{h}{T}\] 
   are homotopic to the identity maps. 
\end{theorem} 

\begin{proof}
The diagram is the direct limit of the diagrams
in Theorem~\ref{Ginduct}~(h) with $h=\lim\limits_{\rightarrow}h_k$,
$g=\lim\limits_{\rightarrow} g_k$ and $f=\lim\limits_{\rightarrow} f_k$, 
where 
\(f_k:\namedright{G_k}{}{\Sigma T^{2np^k}}\) 
is a right inverse for $g_k$ given by Theorem~\ref{Ginduct}~(j).
\end{proof} 

\begin{theorem} 
   \label{Welements}
   The following space belong to $\SW^{\infty}_{r}$:  
   $\Sigma^2\Omega G$, $\Sigma G$, $G\wedge T$, 
   $\Sigma T\wedge T$, and $W$ where $\Sigma T\simeq G\vee W$.
\end{theorem} 

\begin{proof}
This follows immediately from the results in Theorem~\ref{Ginduct} 
by taking limits.
\end{proof} 

The retraction of $T$ off $\Omega G$ in Theorem~\ref{Glimit} induces 
an $H$-structure on $T$ by the composite 
\[m:\namedddright{T\times T}{\tilde{g}\times\tilde{g}} 
    {\Omega G\times\Omega G}{}{\Omega G}{h}{T}.\] 
The following proposition establishes the $H$-fibration property 
in Theorem~\ref{maintheorem}~(a) as a consequence of a slightly 
stronger result. 

\begin{proposition} 
   \label{TisH}
   The map 
   \(\namedright{T}{E}{\Omega S^{2n+1}\{p^{r}\}}\) 
   is an $H$ map with respect to the $H$-space structure $m$ on $T$. 
   Consequently, there is an $H$-fibration sequence 
   \(\nameddright{S^{2n-1}}{}{T}{}{\Omega S^{2n+1}}\). 
\end{proposition} 

\begin{proof}
Filling in the fibration diagram in Theorem~\ref{Glimit} on the right, 
we obtain a homotopy commutative square 
\[\diagram 
     \Omega G\rto^-{h}\ddouble & T\dto^{E} \\ 
     \Omega G\rto^-{\Omega\varphi} & \Omega S^{2n+1}\{p^{r}\}. 
  \enddiagram\] 
Now consider the following diagram 
\[\diagram 
     T\times T\rto^-{\tilde{g}\times\tilde{g}}\drto_{E\times E} 
       & \Omega G\times\Omega G\rto\dto^{\Omega\varphi\times\Omega\varphi} 
       & \Omega G\rto^-{h}\dto^{\Omega\varphi} & T\dlto_{E} \\ 
     & \Omega S^{2n+1}\{p^{r}\}\times\Omega S^{2n+1}\{p^{r}\}\rto 
       & \Omega S^{2n+1}\{p^{r}\}. & 
  \enddiagram\] 
The middle square commutes as $\Omega\varphi$ is an $H$-map and we 
have just seen that the right triangle commutes. The left triangle 
commutes since $\varphi \sim Eh$, so $\varphi\tilde{g}\sim E$. As 
the top row is the definition of the multiplication $m$ on $T$, the 
commutativity of the diagram implies that $E$ is an $H$-map. 

Consequently, the composition 
\(\nameddright{T}{E}{\Omega S^{2n+1}\{p^{r}\}}{}{\Omega S^{2n+1}}\) 
is an $H$-map as it is a composite of $H$-maps, and so the homotopy 
fibration 
\(\nameddright{S^{2n-1}}{}{T}{}{\Omega S^{2n+1}}\) 
is of $H$-spaces and $H$-maps.   
\end{proof} 

The next proposition and the following corollary give structural 
properties of the spaces $T$, $G$, and $R$. 

\begin{proposition} 
   \label{GTatomic}
   The spaces $T$ and $G$ are atomic.
\end{proposition} 

\begin{proof}
It is easy to see that $T$ is atomic using the product
structure and the Bockstein relations. The case of $G$ 
is more difficult. We first show that if $G$ is not atomic then the map 
\(\namedright{P^{2np^k}(p^{r+k})}{\alpha_{k}}{G_{k-1}}\) 
is null homotopic for some $k$. Suppose
\(\gamma:\namedright{G}{}{G}\) 
is a map with the property that
\(\gamma\vert_{G_{k-1}}:\namedright{G_{k-1}}{}{G_{k-1}}\) 
is a homotopy equivalence and $\alpha_k$ has
order $p$. Consider the diagram 
\[\diagram 
     P^{2np^k}(p^{r+k})\rto^-{\alpha_{k}}\dto^{d} 
         & G_{k-1}\rto\dto^{\gamma} & G_{k}\dto^{\gamma} \\ 
     P^{2np^k}(p^{r+k})\rto^-{\alpha_{k}} & G_{k-1}\rto & G_{k}. 
  \enddiagram\]
Since 
\(\gamma\vert_{G_{k-1}}:\namedright{G_{k-1}}{}{G_{k-1}}\) 
is an equivalence, $\gamma \alpha_k$ has the
same order as $\alpha_k$. Consequently $d\not\equiv 0$ (\mbox{mod $p$}) 
and hence $d$ is an equivalence. It follows that 
\(\gamma\vert_{G_k}:\namedright{G_k}{}{G_k}\) 
is an equivalence.

Suppose now that $\alpha_k\sim \ast$. Then we can construct a map 
\[s:\namedright{P^{2np^k+1}(p^{r+k})}{}{G_k}\]
which induces an isomorphism in $H^{2np^k+1}(\ )$. We now show
that the composite 
\[\nameddright{H^{2np^k}(T)}{h_k^{\ast}}{H^{2np^k}(\Omega G_k)} 
    {(\Omega s)^{\ast}}{H^{2np^k}(\Omega P^{2np^k+1}(p^{r+k}))}\]
is an isomorphism. First we note that the image 
of $(h_k)^{\ast}$ is a direct summand since 
\(g_{k+1}:\namedright{T^{2np^{k+1}-2}}{}{\Omega G_k}\) 
induces a left inverse in cohomology. However, there is an isomorphism 
\[H^{2np^k}(\Omega G_k)\cong 
   H^{2np^k}(\Omega G_{k-1})\oplus H^{2np^k}(\Omega P^{2np^k+1}(p^{r+k}))\]
since the inclusion 
\(\namedright{G_{k-1}\vee P^{2np^k+1}(p^{r+k})}{} 
     {G_{k-1}\times P^{2np^k+1}(p^{r+k})}\) 
is an equivalence in this range. However, $H^{2np^k}(\Omega G_{k-1})$
contains no elements of order $p^{r+k}$ since 
$\Sigma^2\Omega G_{k-1}\in\SW_r^{r+k-1}$
by Theorem~\ref{Ginduct}~(i). We conclude that the composition is an
isomorphism. Now apply cellular approximation to
obtain a homotopy commutative diagram 
\[\diagram 
     P^{2np^{k}}(p^{r+k})\rrto^{s^{\prime}}\dto^{E} & & T^{2np^{k}}\dto \\ 
     \Omega P^{2np^k+1}(p^{r+k})\rto^-{\Omega s} & \Omega G_{k}\rto & T  
  \enddiagram\] 
where $E$ is the suspension and $s^{\prime}$ is a skeletal factorization. 
It follows that $s^{\prime}$ induces an isomorphism in
$H^{2np^k}(\ )=\Z/p^{r+k}$.
From this we see that
\[\namedright{H_{2np^k-1}(P^{2np^k}(p^{r+k}):\Z/p)}{(s^{\prime})_{\ast}} 
    {H_{2np^k-1}(T;\Z/p)}\] 
induces an isomorphism as well because of the Bockstein
structure. However $H_{2np^k-1}(T; \Z/p)$ is generated by $u_k$
which is decomposable if $k>0$. This is a contradiction
which implies that $\alpha_k$ is essential and $G$ is atomic.
%page 34
The fact that $G$ and $T$ are atomic and the maps $f$, $g$, $h$
exist implies that $(G,T)$ is a corresponding pair in the sense
of~\cite{G7}.
\end{proof} 

\begin{corollary} 
   \label{Rwedge}
   $R\in \SW_r^{\infty}$. 
\end{corollary} 

\begin{proof}
According to \cite[Theorem~3.2]{G7}, $R$ is a retract of
$\Sigma T\wedge T\in \SW_r^{\infty}$.
\end{proof} 

The next proposition implies that the space $T$ constructed in this 
paper is homotopy equivalent to the space Anick constructed in~\cite{A} 
when $p\geq 5$ (the primes for which Anick's construction holds).  

\begin{proposition} 
   \label{Tunique}
   Suppose $X$ is an $H$ space and there 
   is a fibration sequence:
   \[\nameddright{\Omega^2S^{2n+1}}{\varphi}{S^{2n-1}}{i}{X}\]
   such that the composite 
   \[\nameddright{\Omega^2 S^{2n+1}}{\varphi}{S^{2n-1}}{E^2} 
       {\Omega^2 S^{2n+1}}\]
   is homotopic to the $p^r$ power map. Then $X\simeq T$.
\end{proposition} 

\begin{proof}
Consider the diagram of fibrations 
\[\diagram 
     \Omega W_{n}\rto\dto & \Omega X\rto\dto 
        & \Omega^{2} S^{2n+1}\{p^{r}\}\dto \\ 
     PW_{n}\rto\dto & \Omega^{2} S^{2n+1}\rdouble\dto^{\varphi} 
        & \Omega^{2} S^{2n+1}\dto^{p^{r}} \\ 
     W_{n}\rto & S^{2n-1}\rto^-{E^{2}}\dto & \Omega^{2} S^{2n+1} \\ 
     & X. & 
  \enddiagram\] 
Since $p\cdot\pi_{\ast}(W_n)=0$ and
$p^r\cdot\pi_{\ast}\left(S^{2n+1}\{p^r\}\right)=0$ we conclude that 
$p^{r+1}\cdot\pi_{\ast}(X)=0$. Since $\pi_{2np-1}(W_n)=0$, we also
see that $p^r\cdot\pi_{2np-1}(X)=0$. According to \cite[Corollary~4.2]{AG}
this is sufficient to construct a map 
\[\varphi:\namedright{G}{}{\Sigma X}\] 
which induces an isomorphism in $\pi_{2n}$. The construction
given in \cite{AG} depends only on the co-$H$ space structure
on $G$ and the fact that $\alpha_k$ is divisible by $p^{r+k-1}$, so
the proof works in this context as well. From this we
construct the composition 
\[\namedddright{T}{g}{\Omega G}{\Omega \varphi}{\Omega\Sigma X}{}{X}.\] 
It is an easy calculation with the Serre spectral sequence that 
$H^{\ast}(X;\Z/p)\cong H^{\ast}(T; \Z/p)$, so this map is a homotopy
equivalence.
\end{proof}

\end{document}